\pdfoutput=1
\documentclass[journal]{IEEEtran}
\usepackage{cite}

\ifCLASSINFOpdf
   \usepackage[pdftex]{graphicx}

\else
   \usepackage[dvips]{graphicx}
   \usepackage{epsfig}

\fi
  \graphicspath{{Images/}}
\usepackage[cmex10]{amsmath}
\usepackage{amssymb}  
\ifCLASSOPTIONcompsoc
  \usepackage[caption=false,font=normalsize,labelfont=sf,textfont=sf]{subfig}
\else
  \usepackage[caption=false,font=footnotesize]{subfig}
\fi
 \usepackage{dblfloatfix}
\begin{document}
%
\title{Multi-Agent Orbit Design For Perception Enhancement Purpose}
%
%
%

\author{Hamidreza~Nourzadeh ~\IEEEmembership{Member,~IEEE,}
        and~John E. McInroy,~\IEEEmembership{Senior~Member,~IEEE}
\thanks{H. Nourzadeh is with the Department
of Electrical, Computer and System Engineering, Rensselaer Polytechnic Institute, Troy,
NY, 12180 USA e-mail:nourzh@rpi.edu.}
\thanks{John~E.~McInroy is  with the Department
of Electrical and Computer Engineering, University of Wyoming, Laramie,
WY, 82071 USA e-mail:mcinroy@uwyo.edu.}
}

\maketitle

\begin{abstract}
This paper develops a robust optimization based method to design orbits on which the sensory perception of the desired physical quantities are maximized. It also demonstrates how to incorporate various constraints imposed by many spacecraft missions such as collision avoidance, co-orbital configuration, altitude and frozen orbit constraints along with Sun-Synchronous orbit. The paper specifically investigates designing orbits for constrained visual sensor planning applications as the case study. For this purpose, the key elements to form an image in such vision systems are considered and effective factors are taken into account to define a metric for perception quality.
The simulation results confirm the effectiveness of the proposed method for several scenarios on low and medium Earth orbits as well as a challenging Space-Based Space Surveillance program application.

\end{abstract}


%
\IEEEpeerreviewmaketitle

\section*{Nomenclature}
\addcontentsline{toc}{section}{Nomenclature}
\begin{IEEEdescription}[\IEEEsetlabelwidth{$q_{resolving}$}]
\item[$D$]camera's lens aperture diameter(m).
\item[$R_{ik}$] $i^{th}$ object's position at sample time $k$.
\item[$O_{jk}$] $j^{th}$ observer's position at sample time $k$.
\item[$V_{ijk}$]  a unit vector pointing from $i^{th}$ object to $j^{th}$ observer at sample time $k$.
\item[$n_{ilk}$] A vector indicating $l^{th}$ outward facing unit normal on object $i$ at sample time $k$.
\item[$\alpha _{lit}^{ijk}$] angle between the viewing direction (from $R_{ik}$ to $O_{jk}$) and sunlight direction vector in $radians$.
\item[$\delta_{ijk}$] distance between object $i$ and observer $j$ at sample time $k$ (m).
\item[$\lambda$] Wavelength of the light (m).
\item[$q_{resolving}$] resolving quality of the lens.
\item[$q_{los}$] effect of occlusion on image quality.
\item[$q_{lum}$] illumination effect on image quality.
\item[$q_{angle}$] light angle effect on the image quality.
\item[$q_{view}$] sides observation quality.
\item[$q_{ijkl}$]  Observation quality metric of observing object $i$ by observer $j$ during sample $k$ along face $l$.
\item[$Q$] 4-D quality array.
\item[$J_{sum}$] Cumulative quality matrix.
\item[$a$] Orbit's semi-major axis.
\item[$e$] Orbit's eccentricity parameter.
\item[$\bar{c}$] Altitude's upper bound.
\item[$\underline{c}$] Altitude's lower bound.
\item[$\iota$] Orbit's inclination parameter.
\item[$\Omega$] Longitude of the ascending node.
\item[$\dot{\Omega}$] Nodal precession rate.
\item[$\omega$] Argument of Perigee.
\item[$\nu$] True anomaly.
\item[$\mu$] Earth gravitational constant.
\item[${}^j X_S$ ]  State vector of $j^{th}$ agent.
\item[${}^i X_R$] State vector of $i^{th}$ target.
\item[$\Phi_l$] Perception function of $l^{th}$ sensible quality.
\item[$p$] Optimization parameter vector.
\item[$r$] Orbiting object position with respect to the Earth center.
\item[$\bar{J_2}$] Dimensionless zonal harmonic coefficient of the Earth.
\item[$J_2$] Second gravity coefficient of the Earth.
\item[$J_3$] Third gravity coefficient of the Earth.
\end{IEEEdescription}

\section{Introduction}
%
%
%
%

\IEEEPARstart{R}ecent developments in technology have led to employing robotic systems with more sophisticated data-gathering equipment consisting of various types of sensors. This enables scientists to utilize agents that can simultaneously perform different tasks on different targets which ultimately improves the system performance, and results in reduction of the mission cost and complexity.

In many practical applications, there is an increasing tendency of using a network of inexpensive agents to increase the system's capability and achieve mission objectives. This excess of sensors develops the ability of exploiting the synergistic operation of the sensors to  achieve more accurate measurement and also makes the mission more robust to possible errors and faults. 
In such a multi-agent system, in order to make use of resources more effectively, a meticulous  plan of action is required to distribute tasks among the agents \cite{SensorNetVladimirova2}.

In this paper,  a special class of multi-agent systems  to acquire sensory measurement from physical phenomena of interest is considered. Many orbit design problems belong to this class of multi-agent systems. In fact, one scenario that frequently happens in various space missions is designing satellites trajectories on which the perception quality of certain physical quantities of some targets are maximized. In this constrained path planning problem, the orbit should also fulfil particular requirements imposed by the mission. Depending upon the system configuration and requirements, the design problem could be quite elaborate even for a small number of the agents and targets.

There is an extensive literature on the orbit design subject for different science missions specifically oriented toward special types of orbits \cite{ABC_SSO,orbitDesignGEOSYNC,OrbitDesignNonKeplerian}. In these studies, the main aim is to define attributes of the mission usually without providing information about multi-agent path planning.  In another interesting field of study, researchers have addressed coordinated trajectory design for multi-agent multi-target systems in the context of formation flying. The formation flying research field seeks to develop strategies to control formation of a satellite fleet in order to meet a global performance objective, and the focus of interest is to employ an adequate control strategy to maintain  the desired formation. \cite{ConvexOrbitDesignMit,Inalhan02relativedynamics,FormationFlying}. 

In many applications, the ultimate goal is to effectively allocate each sensor to a target at any point in time \cite{Hamid2011,HamidRobust2013}. But, before performing the resource allocation process, there are fundamental questions that need to be addressed early in the planning stage. The first and  most important one is how the agents (leaders in cluster formation) should be placed or moved with respect to the targets during the operation to optimize that particular mission \cite{HamidPatrolling2013}. 

In this article, the multi-agent orbit design procedure is formalized as a constrained parametric optimization problem by taking into account different constraints that might be involved in a science mission. Particularly, the SSO and Frozen orbit constraints that are widely employed in practical application are defined in term of optimization parameters.

To analyze the performance of the proposed method, a special Spaced-Based Space Situational Awareness application is chosen. In this application, the main objective is to agents trajectories on which the observation quality of some Resident Space Objects (RSOs) is maximized.
Although the orbit design procedure is performed in the context of inspection quality, the proposed formulation is quite general, and any type of physical quantity can be incorporated and optimized.
To handle possible perturbations and uncertainties in the system, a max-min model is employed to bring robustness into the planning, and, to deal with large-scale cases,   some guidelines are proposed to choose an adequate numerical optimization algorithm.

This paper is organized as follows. Section \ref{sec:Formulation} formulates the  coordinated trajectory  design problem that determine paths on which the quality of received information is maximized in a sense.  Section \ref{sec:QualityCalc} discusses contributing elements in image formation, and defines a quality metric for a Space Situational Awareness (SSA) application by considering the dominating factors. Section  \ref{sec:OrbitalElements} describes the motion equations of the objects and  explains two sets of parameters commonly used to uniquely identify an orbit.  In Section \ref{sec:ProblemConstraints}, the constrained multi-agent multi-target trajectory optimization problem is constructed, and it is described how to incorporate the constraints of different configurations that might be employed in a science mission.  Brief remarks on choosing suitable numerical algorithms to solve the optimization problem are given in Section \ref{sec:NumericalOptAlg}.  Simulation studies are included in Section \ref{sec:Sim}, and the conclusions are drawn in the final section.

\section{Problem Formulation}
\label{sec:Formulation}

\begin{figure}[!htb]
\includegraphics[width=.35\textwidth]{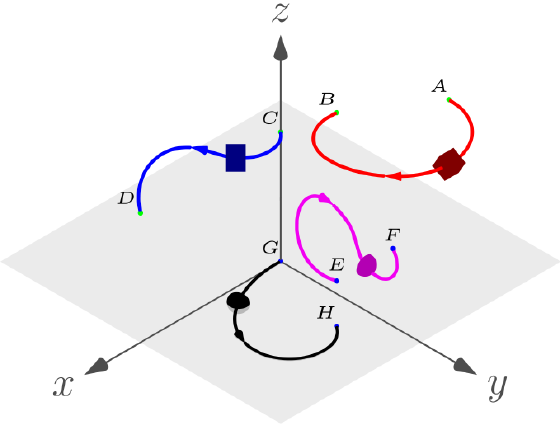}
\centering
\caption[illustrative constraint path planning example]{Path planning to maximize sensory reception.}
\label{fg:PathPlanning}
\end{figure}

This section focuses on the problem of designing the best agent trajectories on which sensory perception is maximized subject to some constraints imposed by the system nature.
Consider a general path planning problem illustrated in Figure \ref{fg:PathPlanning}. As depicted, two cube-shaped objects are moving on two different curved paths (AB and CD). The objects passively or actively send out certain types of signals that can be perceived at any point of the space with different quality levels. At a particular point in time, depending on the system configuration, specific quality fields are formed for any perceptive physical aspects of the objects. In fact, when objects move on their paths, a time dependent vector or scalar field is formed for any of these physical quantities.

In the case of a scalar field, given the state of target $i$(${}^i X_R$), state of sensor $j$(${}^j X_S$) and effective physical factors, the total achievable quality along the sensor trajectories, e.g. EF and GH, in the space is the sum of values of the perceptions at all points on the curve, and it can be expressed by Equation \ref{eq:lineIntegral1}.


\begin{IEEEeqnarray}{rCl}
\label{eq:lineIntegral1}
& & \quad i=1,\cdots,n \nonumber\\*[-0.45\normalbaselineskip]
\int\limits_{t_1}^{t_2}{\Phi_l({}^i X_R(t),{}^j X_S(t)) dt}, & & \quad  j=1,\cdots,m
\\*[-0.45\normalbaselineskip]
&  & \quad l=1,\cdots,L \nonumber
\end{IEEEeqnarray}

Where $\Phi_l$ is a function that describes the perception value of quantity $l$ sensed by sensor $j$ from the object $i$ at any point of time. The symbols $m$,$n$ and  $L$ denote the number of agents, targets and physical quantities, respectively.

In the discrete case, the perception quality of target $i$ obtained by all sensors can be written as Equation (\ref{eq:lineIntegral2}). For $i=1,\cdots,n$ and $l=1,\cdots,L$,  (\ref{eq:lineIntegral2}) gives a perception quality matrix $J_{sum} \in \mathbb{R}^{n \times L}$.
The $(i,l)$ element of $J_{sum}$ gives the $l^{th}$ quality attribute of the $i_{th}$ object along the trajectories.
\begin{equation}
\label{eq:lineIntegral2}
J_{sum}=\sum\limits_{j=1}^{m}{\sum\limits_{k=1}^{N}{\Phi_l({}^i X_R[k],{}^j X_S[k])}}.
\end{equation}

In this paper, it is presumed that each trajectory can be parametrized with a set of parameters. Concretely, $^{j}X_S$ moves on the trajectory $j$ that is parametrized by $p_j  \in \mathbb{R}^{np_{j}}$ where $np_{j}$ is the number of parameters required to represent trajectory $j$. Let $p \in \mathbb{R}^ {\sum\limits^{m}_{j=1} np_j} $ be a vector consisting of all the requisite parameters to describe the sensors trajectories (Equation (\ref{eq:ConcatOptParam1})).

 \begin{equation}
\label{eq:ConcatOptParam1}
	{{p}}=\underset{j=1}{\overset{m}{\mathop{\oplus }}} p_j.
 \end{equation}

Here, $\oplus$ is the vertical concatenation symbol.

It is worthwhile noting that the relation between parameters $p_j$  and trajectory $j$ doesn't have to be explicit, and even a heuristic path planning algorithm that generates trajectory $j$ based on tuning parameter $p_j$ can be utilized here.

Depending upon the application and how to maximize the $J_{sum}$ matrix, various possible optimization strategies can be considered to acquire different sets of sensors' trajectories. In this article, to deal with worst cases, a max-min optimization model is employed to design trajectories which enforce  planning redundancy, and make the result robust to the probable uncertainty sources in the system. Therefore, the optimization formulation determines the sensors' trajectories that maximize the minimum achievable perception for all objects (the minimum element of $J_{sum}$) subject to some constraints imposed by the system (Equation (\ref{eq:MaxMin1})). 
Figure \ref{fg:PathPlanning} portrays the case when each agent's trajectory is confined to planes $z=0$ and $x=0$.

\begin{equation}
\label{eq:MaxMin1}
  \left\{
  \begin{array}{l l}
   \underset{p}{ \text{Maximize}} &   min(J_{sum})\\
    \text{subject to} &   : \text{Constraints}
    \end{array} \right.
 \end{equation}

As mentioned above, two prime requirements to construct the optimization problem  \ref{eq:MaxMin1} are the perception function, $\Phi$, and the trajectories parameter vector, $p$. In two subsequent sections, these elements are defined for a visual measurement SSA application. In this scenario, it is assumed that all the agents are equipped with passive visual sensors and the physical quantities that are maximized along the trajectories are inspection qualities of the different sides of the objects. Furthermore, each trajectory is parametrized by orbital elements parameters. While the following section elaborates the calculation process of the observation quality given states of the objects and the cameras in this scenario, Section \ref{sec:OrbitalElements} describes motion dynamics of agents and targets as well as  the trajectories parametrization process.

\section{Optical Imaging Quality Calculation}
\label{sec:QualityCalc}
There are many vision systems with different sensors and ways of data acquisition, and the process of sensor planning highly depends on the way that information is gathered in that specific application.
In some vision systems, only passive sensors are used while both passive and active  sensors are used in some other machine vision systems. For simplicity, suppose that all vision systems in our application are equipped with passive sensors although the methods developed here apply also active sensing. To compute the qualities of the resulting images, all effective factors for forming an image should be considered, i.e. the intrinsic parameters of the sensors and environmental factors. There are a lot of techniques for automatically determining the image quality \cite{irv97a,irv03a,lea96a,lea97a,mav95a}.

In this section, the preliminary requirements to compute the observation quality in an optical imaging system will be demonstrated in detail by investigating all effective aspects in forming images and defining a quality metric based upon them.

\subsection{Effective Factors on Image Quality }
\begin{figure}[!htb]
\centering
\includegraphics[width=.5\textwidth]{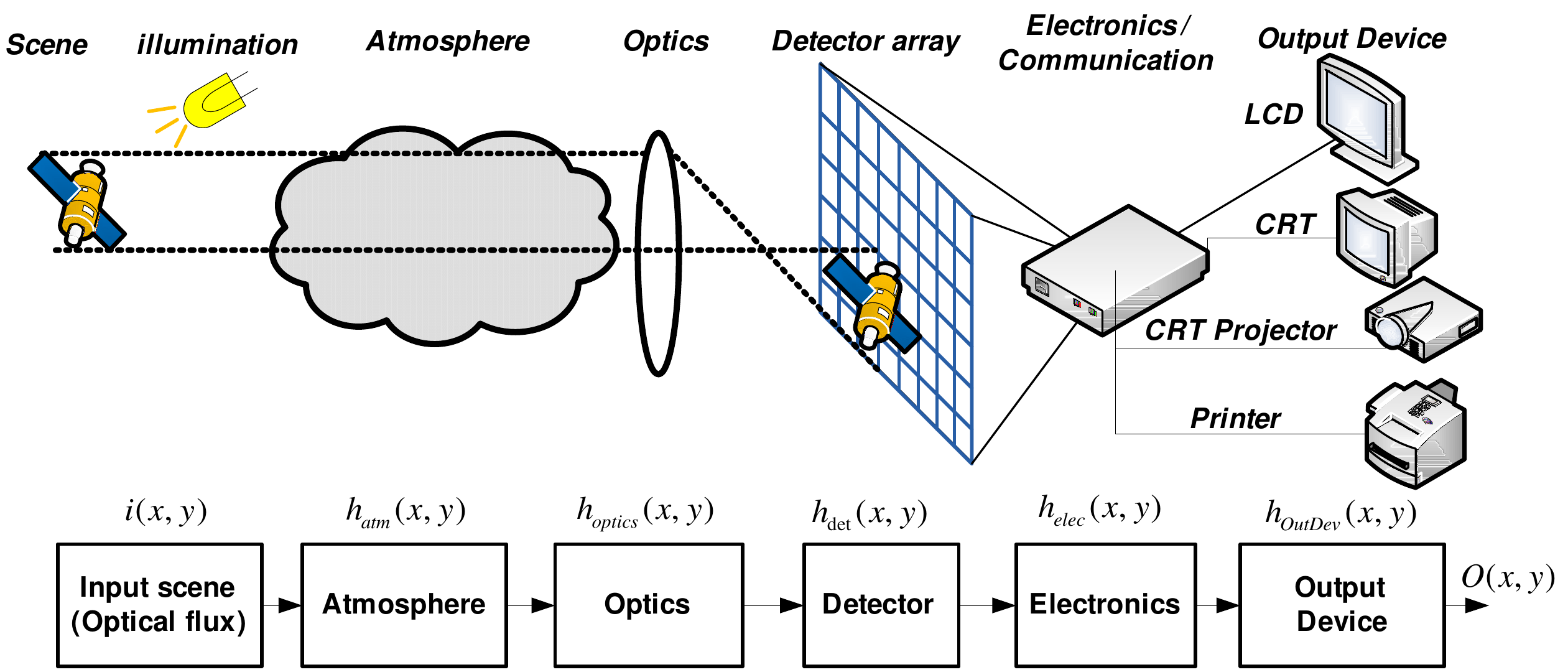}
\caption[Imaging System Components]{A comparison of the modulation transfer
function and the contrast transfer function.} \label{fg:ImagingSysComp}%
\end{figure}

The output of an optical imaging system depends on several physical factors, i.e. lighting, atmospheric attenuation, light diffraction, occlusion, imaging sensor resolution and sensitivity, electronics parts and output devices (Figure \ref{fg:ImagingSysComp}).  In the case of a linear shift-invariant imaging system, the output image is mathematically expressed by the convolution of the Point Spread Function (PSF which is the impulse response of the imaging system) and input. Moreover a typical imaging system can be shown as impulse response blocks of some physical effects in series where the input signal is optical flux from scene and the output is the resulting image \cite{young_signal_2008}. Figure \ref{fg:ImagingSysComp} depicts a typical imaging system.

The Fourier transform of PSF is called  the \textit{optical transfer function} (OTF), and it is the product of the component OTFs (Equation (\ref{eq:OTFImagingSys})).


\begin{IEEEeqnarray}{rCl}
\label{eq:OTFImagingSys}
O_{sys}(f_x,f_y) &=& I(f_x,f_y) H(f_x,f_y) \nonumber\\*[-0.1\normalbaselineskip]
H(f_x,f_y)&=& H_{atm}(f_x,f_y)H_{opt}(f_x,f_y) H_{det}(f_x,f_y)\\*[-0.1\normalbaselineskip]
&  &   H_{elec}(f_x,f_y)H_{dev}(f_x,f_y) \nonumber
\end{IEEEeqnarray}

In Fourier optics, The Modulation Transfer Function (MTF) is defined as the magnitude of the OTF. MTF is extensively used to assess the resolving ability of an imaging system, and it is the most prevalent scientific method to predict the quality of an imaging system. To perform empirical MTF measurements, a test target which contains black/white line pairs is commonly used. The test target corresponds to a square wave rather than a sine wave which is used in MTF computation. This empirical method is called  the Contrast Transfer Function (CTF) which overestimates MTF, and it can be utilized to evaluate the performance of an imaging system \cite{ModulationTransferFunction}.

In SSA applications,  only some physical aspects have a significant influence on the output image. So, in  subsequent sections, only illumination, light fraction, occlusion and side observation quality will be taken into account.

\subsubsection{Light diffraction}

In optics, diffraction is categorized into two different classes, i.e. Fraunhofer  diffraction, or far-field diffraction, and Fresnel diffraction, or near-field diffraction. Fraunhofer diffraction occurs when waves from far-field distance are passed through an aperture or slit. This happens when Fresnel number $F \ll 1$  and the parallel rays approximation is applicable.

Because of the Fraunhofer diffraction phenomenon, even an optical lens with perfect quality has limited performance, and the MTF of each lens expresses this limitation in spatial frequency.  In fact, light coming from a point light source diffracts and forms a Airy disk pattern. Using this pattern,  Rayleigh  proposed a criterion to find a optical resolution of the lens. The criterion states that two points with a angular separation equal to the angular radius of the Airy disk can be resolved. Equation \ref{eq:Rayleigh} relates the angular resolution $ \theta$ to wave length of light ($\lambda$) and diameter of lens aperture $D$ \cite{hecht2002optics}.

\begin{eqnarray} \label{eq:Rayleigh}
sin \theta = 1.22 \frac{ \lambda}{D}
\end{eqnarray}

 The image resolving quality for a system with a fixed aperture diameter $D$ and fixed wavelength $\lambda$ can be expressed by Equation (\ref{eq:QualityBlure}), where $\delta$ is the distance between the object and the lens.

\begin{eqnarray} \label{eq:QualityBlure}
q_{resolve}= \frac{D}{1.22 \lambda  \delta}
\end{eqnarray}

Thus, by increasing the object distance, the resolving ability decreases because of the diffraction.

\subsubsection{Illumination}

The relationship between reflected light (radiance) and incoming illumination depends on the direction in which light arrives, as well as the shape and type of surface. Absorption, transmission, scattering or a combination of these effects occurs to the incoming light when it strikes the surface. In other words, the intensity and colour reflected depends on the illumination and reflection angle and the surface material. Surfaces can be categorized into three different groups, i.e. diffuse, Lambertian and specular surfaces. Specular surfaces behave like a mirror and reflect light into a lobe of the specular direction  so the reflected light highly depend on  illumination direction, while for Lambertian and diffuse surfaces such as cotton cloth, matte paper and matte paint the radiance leaving the surface does not have any meaningful correlation with illumination direction, and the Bidirectional Reflectance Distribution Function (BRDF) of these surfaces is constant \cite{forsyth_computer_2002}.

\begin{figure}[!htb]
\centering
\includegraphics[width=.35\textwidth]{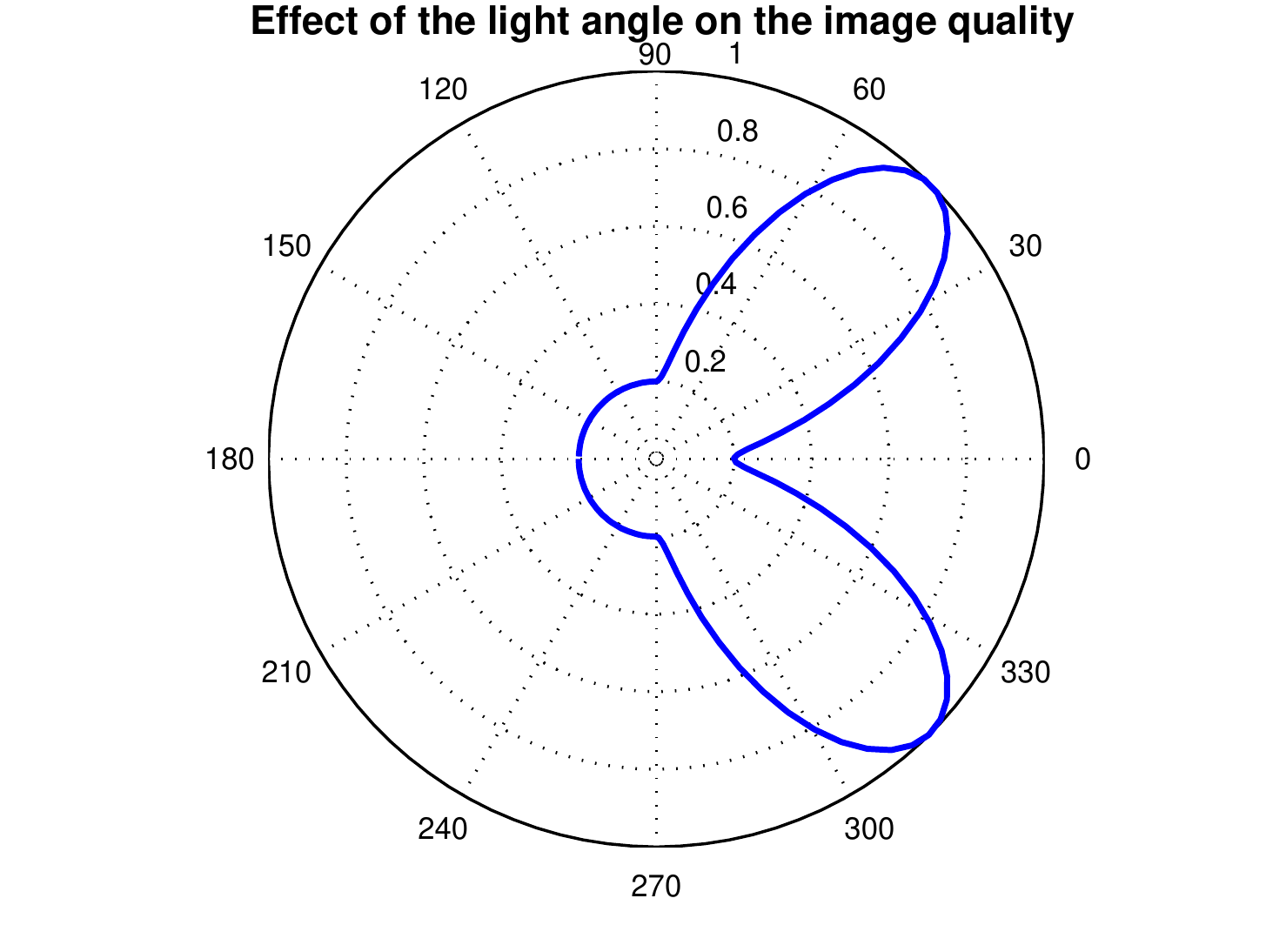}
\caption[Illumination Quality ]{Illumination quality of a specular surface for different values of observation angle .} \label{fg:Qlit}%
\end{figure}

\begin{figure}[!htb]
\centering
\includegraphics[width=0.4\textwidth]{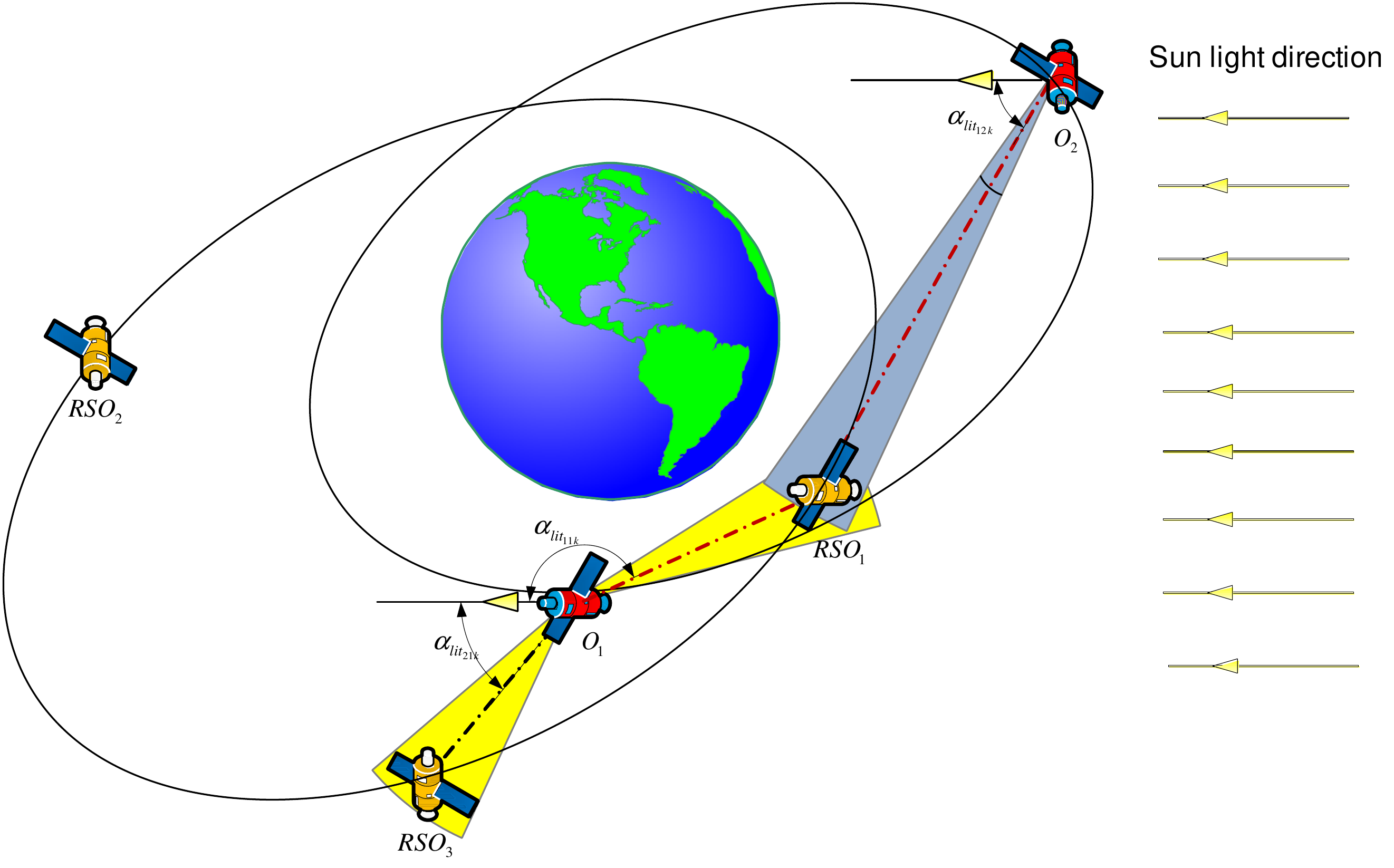}
\caption[Finding illumination quality in a SSA application]{Finding illumination quality in a SSA application.} \label{fg:satIllum}%
\end{figure}

Since the outer surfaces of man-made satellites are highly reflective, the angle of the light with respect to the viewing axis is extremely critical. Let $\alpha_{lit}$ be the angle between the viewing axis
and incoming sunlight (Figure \ref{fg:satIllum}). When $\alpha_{lit} \simeq 0 $, the image contains much specular glare. On the other hand, if $\alpha_{lit} > 90^{\circ} $, the satellite will be back-lit by the Sun, and means that it will produce a silhouette with less details. $\alpha = 45^{\circ} $ will provide the best illumination level.
A suitable function that describes the illumination effect and lies between 0 and 1 can be written as Equation (\ref{eq:lightAngle}). Where  $\underline{q_{lit}}$ is the minimum value that illumination quality takes.   Figure \ref{fg:Qlit} depicts values $q_{lit}$  take on  as $ \alpha _{lit} $ changes from 0 to 360 degrees in polar coordinate when $\underline{q_{lit}}=0.2$ .

\begin{equation}
\label{eq:lightAngle}
q_{lit}= \left\lbrace
\begin{array}{l c}
\underline{ q_{lit}} + (1-\underline{ q_{lit}}) sin^2(2\alpha_{lit})   &  \;\; \alpha_{lit}<\frac{\pi}{2} \\
\underline{ q_{lit}} & \;\; \alpha_{lit} \geq \frac{\pi}{2}\\
\end{array}
\right.
\end{equation}

Figure \ref{fg:satIllum} illustrates the illumination quality of observation for three Resident Space Objects (RSOs) which can be observed by two observer satellites ($O_1$, $O_2$ ). Observer 2 inspects RSO 1, and the angle $ \alpha^{12k}_{lit} $, angle between sunlight and observation direction is  around  $45^{\circ}$. Therefore the illumination quality will be around 1. Observer 1 has two options to inspect, either RSO 1 or RSO 3. The image taken of RSO 3 has a better quality in terms of illumination ($ \alpha^{21k}_{lit}\thickapprox 45^{\circ}$),  while the  the RSO 1 image would have silhouette problem  since $ \alpha^{11k}_{lit}> 90^{\circ}$.

\subsubsection{Occlusion and Sunlight Quality}

Since an occluded object can't be inspected by the observer, the resulting image quality for that object should be considered zero. On the other hand, if there is no obstacle between the object and observer the possibility of having perfect image of that object exists (Equation \ref{eq:qlos}).

\begin{eqnarray} \label{eq:qlos}
   \left\{
  \begin{array}{l l}
    q_{los}=0 & \quad \text{occlusion}\\
    q_{los}=1 & \quad \text{no occlusion}\\
  \end{array} \right.
\end{eqnarray}

In order to observe the object, an acceptable level of illumination is required.  For instance, in satellite imaging , sunlit images are more desirable, and dark images obtained from shadow areas behind the Earth aren't as useful as sunlit images. To incorporate the sunlight quality factor, let $ q_{lum} = 1$ when the object is in sunlight. Otherwise  $q_{lum}$ is a small positive number(Equation \ref{eq:qsun}).

\begin{eqnarray} \label{eq:qsun}
   \left\{
  \begin{array}{l l}
    q_{lum}=1 & \quad \text{lit image}\\
    q_{lum}\ll 1 & \quad \text{dark image}\\
  \end{array} \right.
\end{eqnarray}

\subsubsection{Sides Observation Quality}
In visual sensor planning, the vast majority of  studies are mainly concerned with finding the vision system parameters  to inspect all sides of an object with a minimum number of observations \cite{chen_active_2008}. One way is to assign an outward facing unit normal for each side of interest, i.e.  let  $ n _{ilk}$ be the $l^{th}$ outward facing unit normal on object $i$ at sample time $k$. To determine if a specific side of an object is viewed, let $ v_{ijk}$ be the unit vector pointing from object $i$  to observer $j$ at sample time $k$. Then the observation quality of  side $l$ of object $i$ can be computed as Equation
(\ref{eq:qsides}).

\begin{eqnarray} \label{eq:qsides}
  q^{ilkl}_{view}= \left\{
  \begin{array}{l l}
     n _{ilk}^T v_{ijk} & \quad    n _{ilk}^T v_{ijk} >0\\
     0 & \quad  n _{ilk}^T v_{ijk} \leq 0\\
  \end{array} \right.
\end{eqnarray}

\begin{figure}[h!]
\centering
\includegraphics[width=.2\textwidth]{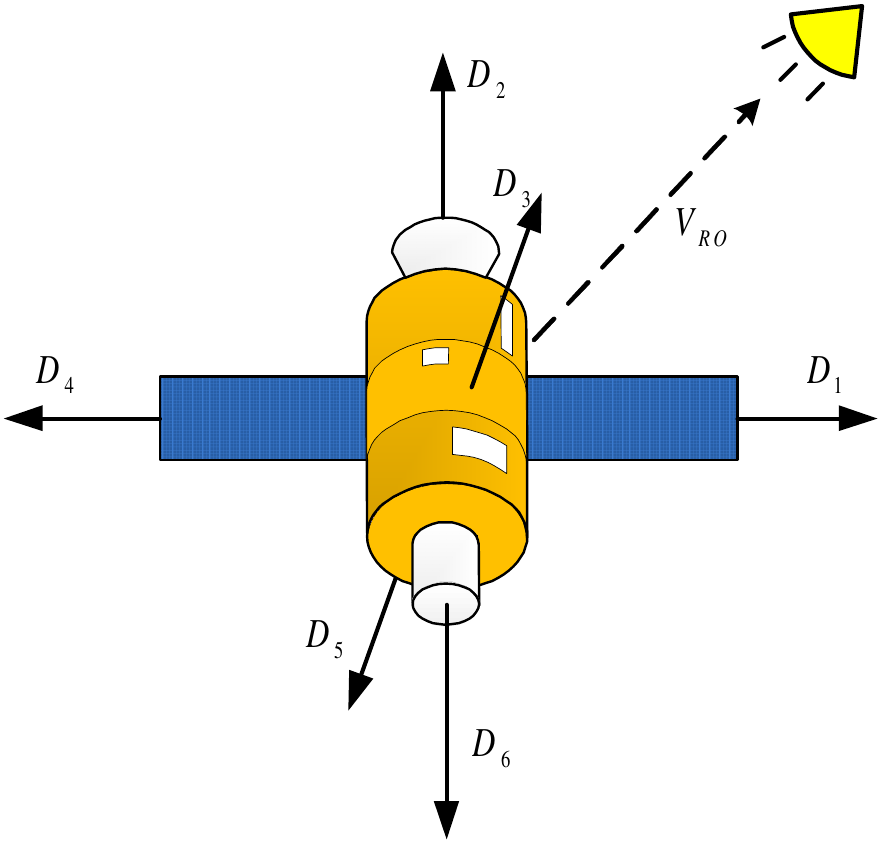}
\caption[Sides observation Quality for a RSO satellite]{Assigned normals to specify sides of interest and sides observation quality computation (inner product).} \label{fg:SidesQuality}%
\end{figure}

Figure \ref{fg:SidesQuality} depicts a RSO with some assigned sides of interest ($D_l$). According to Equation (\ref{eq:qsides}), the ($D_1$) and ($D_2$) have a positive quality and the inner product of $v$ and the other side normals are less and equal to zero. This means that an observer in this direction can partially inspect two sides.

\subsection{ Single Quality Metric }

In order to form a single number reflecting overall observation quality of object $i$ inspected by observer $j$ during sample $k$ along face $l$, all calculated qualities should be combined. Since overall quality varies from application to application, it can be done by multiplying weighted qualities. This method allows us to control the participation of all constituent parts in a single quality metric. Let $ q_{ijkl}$ be the product of all mentioned qualities (Equation (\ref{eq:qtotal})).

\begin{eqnarray} \label{eq:qtotal}
 \Phi_l({}^{i}X_R[k],{}^{j}X_j[k])= q_{ijkl} = q^{ijk}_{resolve}  q^{ijk}_{los}  q^{ijk}_{lum}  q^{ijk}_{angle}  q^{ijkl}_{view}
\end{eqnarray}

For different values of $i$, $j$, $k$ and $l$, $q_{ijkl}$ forms a 4D observation quality array $Q \in \mathbb{R_+} ^{n \times m \times N\times L }$.

\section{Orbital Elements and Motion Equations}
\label{sec:OrbitalElements}

\subsection{Motion Dynamics}
Assuming a uniform spherical Earth,  Equation (\ref{eq:MotionDynamics}) represents the motion equation of the target/object on its orbit.

\begin{equation}
\label{eq:MotionDynamics}
\ddot{\vec{r}}= -\frac{\mu}{r^3} \vec{r} + \vec{d}
\end{equation}

Where $\vec{r} =[x \; y \; z]^T $ is the relative position of the target with respect to the center of the Earth. The vector $\vec{d} = [ d_x \; d_y \;  d_z]^T$ is exogenous forces including  attitude control system (ACS) actions and  modeled or unmodeled  disturbing forces affecting target such as solar pressure, atmospheric drag  etc.  Earth gravitational constant $\mu \triangleq 398600.4418 \;\rm{Km}^3 \rm{s}^{-2}$, and $r=\parallel \vec{r} \parallel _2$ represents the Euclidean distance of vector $\vec{r}$.

Letting state vector $X=[ \vec{r}  \; \; \dot{\vec{r}} ] ^T$ and velocity vector $\vec{v}=\dot{\vec{r}}=[\dot{x} \; \; \dot{y} \; \; \dot{z} ]^T  $, Equation (\ref{eq:MotionDynamics}) can be rewritten as a set of first order ordinary differential equations (ODEs) consisting of only the orbital state vector (Equation (\ref{eq:StateEquation})).

\begin{equation}
\label{eq:StateEquation}
\dot{X}=f(X,\vec{d})=\left[ \begin{matrix}
 \vec{r} \\
 -\mu \frac{\vec{r}}{r^3} +\vec{d} \\
\end{matrix}  \right]
\end{equation}

The discretized motion equation can be numerically solved by choosing an adequate ODE solver (Equation (\ref{eq:DiscreteStateEquation})).

\begin{equation}
\label{eq:DiscreteStateEquation}
X[k+1]=F(X[k],\vec{d}[k]), \; \; X[k] \approxeq X(kT_s).
\end{equation}
Where $T_s$ denotes sampling rate.

\subsection{Orbital Elements}
When the ACS actions and perturbing forces cancel each other, i.e. perfect control, the $\vec{d}$ will be a zero vector, and the target trajectory forms an elliptic orbit. This orbit can be uniquely determined by knowing the orbital elements at an epoch (time $t$). There are several sets of orbital parameters that result in the same motion dynamics.  In classical two body systems, Keplerian orbital elements and  the orbital state vector $X$ are two commonly used sets of orbital elements. The orbital state vector contains two physical quantities, position and velocity, of the moving target at epoch written in the reference frame.

\begin{figure}[!htp]
\begin{center}
\includegraphics[width=.5\textwidth]{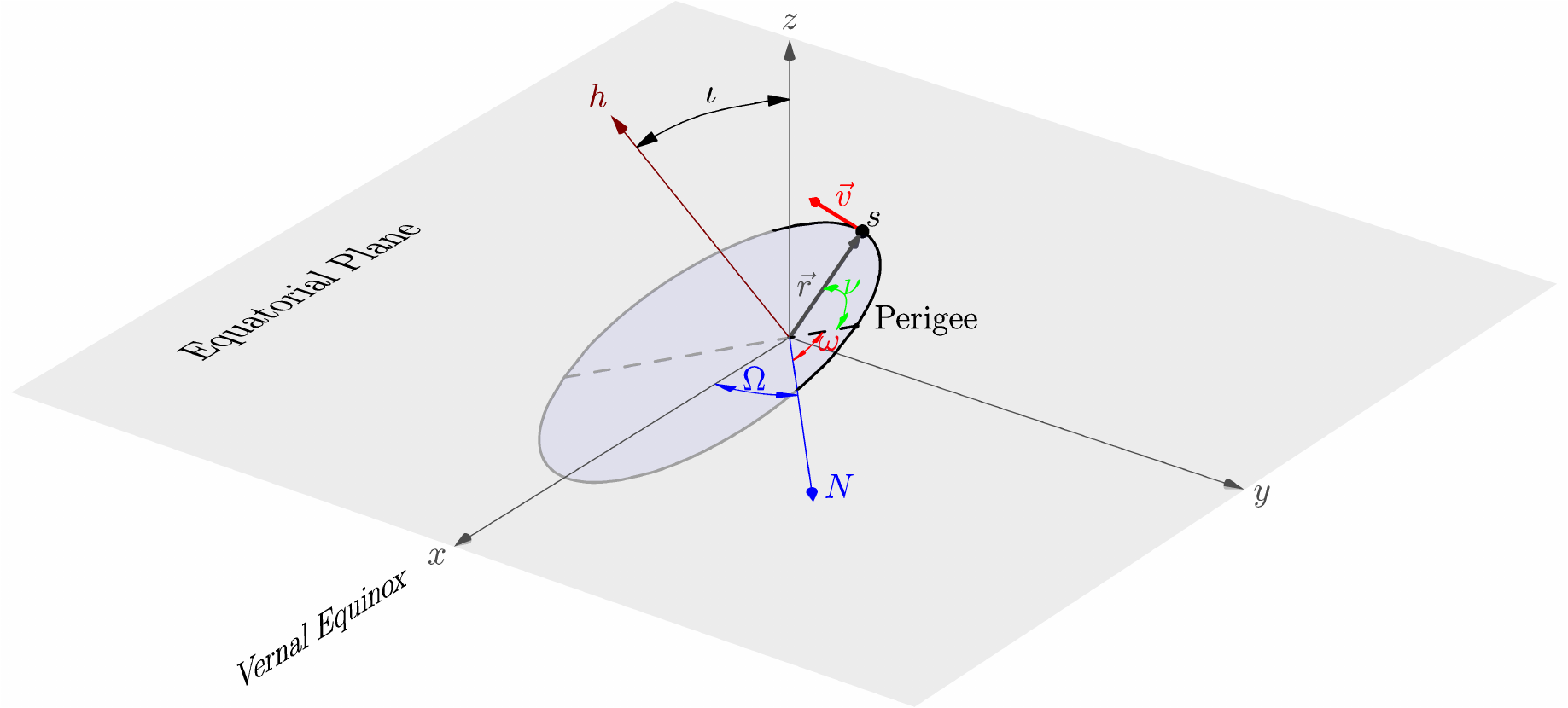}
\caption{Keplerian and Orbital State Elements.
\label{fig:OrbitalElements}}
\end{center}
\end{figure}

On the other hand, the Keplerian orbital elements are directly related to the shape of the resultant orbit and the position of the target on the path. They consist of six parameters: Semi-major axis ($a$) is the sum of the periapsis (perigee) and apoapsis (apogee) distances divided by two, Eccentricity ($e$) controls the shape of the elliptic orbit, and specifies how much the orbit is deviated from a circular orbit, Inclination $\iota$ is the angle between the equatorial plane and the plane containing the elliptical orbit(elliptical plane). Longitude of the ascending node ($\Omega$)  is the angle between ascending node vector and the direction of the Vernal Equinox.  Argument of Perigee ($\omega$) orients the elliptical orbit in the elliptical plane, and it is the angle between the right ascending node direction and the Perigee's direction.
True anomaly ($\nu$) at epoch $t$ represents the geometric angle of the target in the orbital plane (Figure \ref{fig:OrbitalElements}).

Since an orbital state vector can be converted to Keplerian orbital parameters by a non-linear transformation at any epoch \cite{curtis_orbital_2005}, throughout this paper, any of these representations is employed wherever it is more convenient to use.

\section{Optimization Formulation}
\label{sec:ProblemConstraints}
The problem of maximizing the minimum observation quality for all sides of the targets with respect to the sensor's trajectories parameters can be written as Equation (\ref{eq:GeneralOptimization}).

\begin{equation}
\label{eq:GeneralOptimization}
  \left\{
  \begin{array}{l l}
   \underset{p}{ \text{maximize}} &  f(p)\\
    \text{subject to} &   g_i(p) \leq 0, \; i=1,\cdots, n_g \\
  &  h_j(p) = 0, \; j=1, \cdots, n_h
    \end{array} \right.
 \end{equation}

Where  $n_g$ and $n_h$ refer to the number of inequality and equality constraints, respectively.

  \begin{equation}
\label{eq:ConcatOptimizatioinParam}
	{{p}}=\underset{j=1}{\overset{m}{\mathop{\oplus }}}\,\,{{p_j}} = \underset{j=1}{\overset{m}{\mathop{\oplus }}}\,\,{[a_j,e_j,\Omega_j,\omega_j,\iota_j,\nu_j]^T }.
 \end{equation}

Subscript $j$ denotes the number assigned to the agent that moves in orbit with the orbital parameter $p_j$, and $f(p)=min(J_{sum}(p))$.

\subsection{Scaling}
Scaling has a profound impact on the performance of many  optimization problems, and  ignoring it can easily degenerate convergence speed and cause numerical problems.
Beside these facts, analysis of the scaled parameters is  much simpler especially in practical applications.
In Equation (\ref{eq:GeneralOptimization}), vector $p$  elements vary in different ranges.
For instance, while $e_j$ can change between 0 and 1, $a_j$ are typically of the order of 7 to 40 thousands Km.
Usually, scaling process can be performed by dividing each element of vector $p$ by its maximum expected range plus its minimum expected value.
let $\hat{p} \in \mathbb{R}_+^{6m}$ be the scaled parameters, and $\bar{p}$ and $\underline{p}$  denote upper bound and lower bound of vector $p$.
Scaling and unscaling processes can be written as Equation (\ref{eq:Scaling}). In this case, all $\hat{p}$ elements belong to the 0 to 1 interval.
\begin{equation}
\label{eq:Scaling}
\hat{p}=(p-\underline{p}) \oslash (\bar{p}-\underline{p}), \quad p= \hat{p} \otimes (\bar{p}-\underline{p})+\underline{p}
\end{equation}
Where $\oslash$ and $\otimes$ are the element wise division and multiplication operator, respectively.

Now one can readily obtain the scaled optimization problem by substitution of $p$ with $\hat{p} \otimes (\bar{p}-\underline{p})+\underline{p}$ in  Equation (\ref{eq:GeneralOptimization}). After solving the unscaled optimization problem, the scaled solution $\hat{p}$ can be converted into the original orbital element parameters $p$.

\subsection{Co-Orbital Configuration}
In many surveillance missions, some agents are moving on the same orbit. In this case, all orbital elements of these agents are the same except $\nu$  that can be freely chosen. To enforce this constraint for agent $j_1$ and ${j}_{2}$ that share the same orbit, the following linear equality constraint should be added to the optimization problem.
\begin{equation}
\label{eq:SameOrbitConstraint}
A^T p =0
\end{equation}

Where $A^T$ is a $5 \times 6m$ matrix consisting of $m$ horizontally concatenated $5 \times 6$ block matrices. All of these sub-matrices are zero matrix except the $j_1 ^{th}$ and  $j_2 ^{th}$ ones that are replaced by $A_s$ and $-A_s$, respectively.
Mathematically, matrix A can be constructed by the following equations:
\[A=\underset{j=1}{\overset{m}{\mathop{\oplus }}}\,\,\,\left\{ \begin{array}{*{35}{l}}
   {{A}_{s}^T} 			& 	\; \; \; j={{j}_{1}}  \\
   -{{A}_{s}^T} 		& 	\; \; \; j={{j}_{2}}  \\
   {{0}_{6\times 5}} 	& 	\; \; \; \mathrm{otherwise}  \\
\end{array} \right.,   \;   A_s=\left[ I_{5 \times 5} \; , \; 0_{5 \times 1} \right] \].

\subsection{Uniform Distribution}
\begin{figure}[h!]
\centering
\includegraphics[width=.3\textwidth]{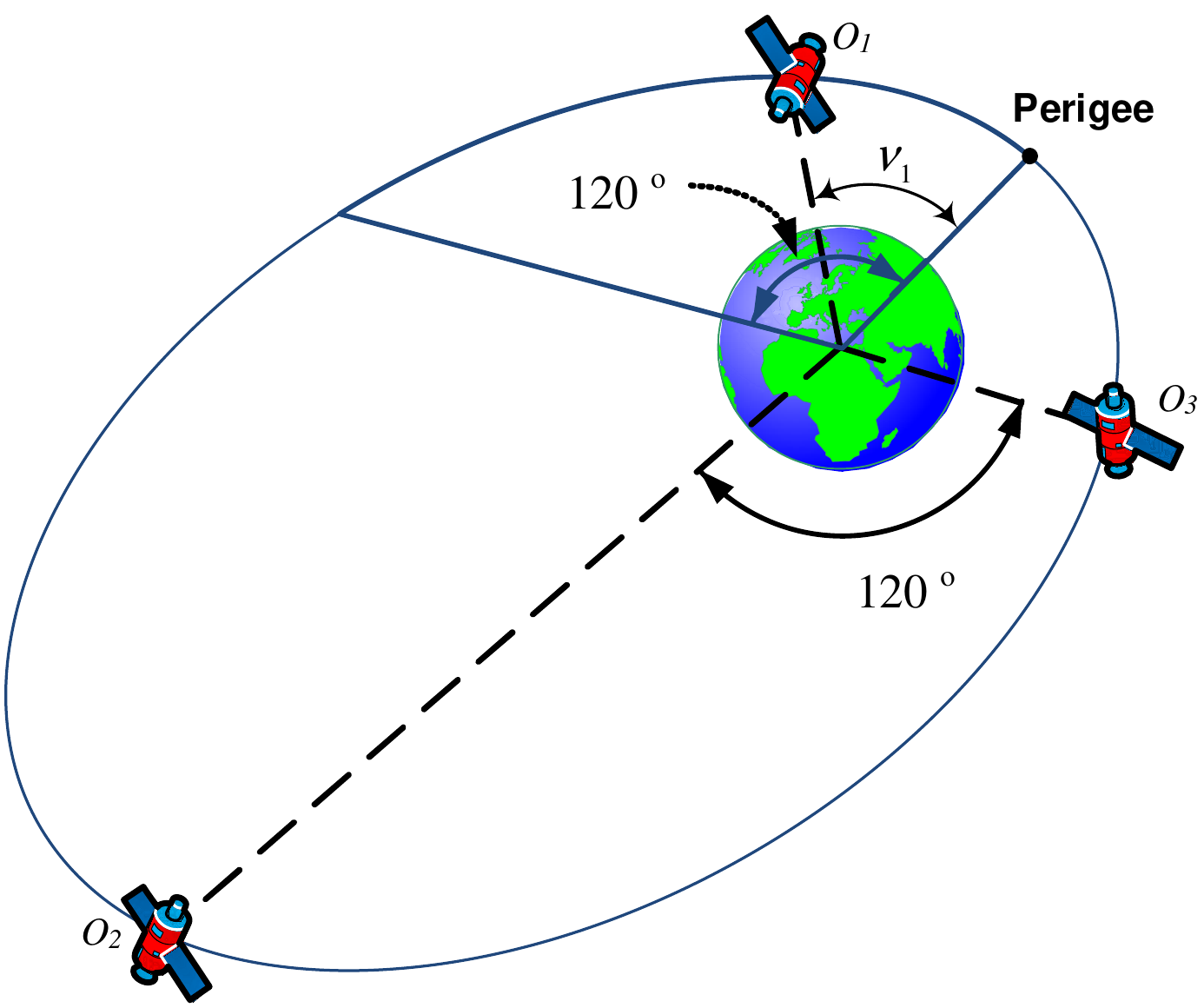}
\caption[Agent distribution in an Orbit]{Sperading satteiltes in an orbit.} \label{fg:DistributeSats}%
\end{figure}

Another configuration that might be considered in orbit design is to uniformly distribute agents which have a common orbit (trailing formation). Let $O$ be the number of distinct observers orbits and  $m_o$ denote the number of satellites in the orbit $o$.  In this case, the true anomaly difference of the initial states of two adjacent satellites in orbit $o$ should be $\frac{360}{m_o}^{\circ}$. It is worthwhile to mention that only in unperturbed circular orbit will the difference be maintained during motion.
Figure \ref{fg:DistributeSats} illustrates one possible configuration when three observer satellites are evenly distributed with a constant true anomaly distances.

To incorporate this scheme in the optimization problem, consider a subset of agents indices that share the same trajectory (orbit $o$) $\{j_1,j_2,\cdots,j_{m_o}\}$.  Assuming different types of satellites with different characteristics,  there will be $(m_o-1)!$ ways to arrange these $m_o$ distinct objects in orbit $o$. To find the optimal arrangement, one should solve the optimization problem for all possible orders and choose the one that has the best resultant cost function. This method is impractical when dealing with a large scale system, and to avoid this cumbersome process, an approximate approach should be taken into account.  One heuristic way that results in satisfactory suboptimal solutions can be obtained by assuming identical agents, and  solving the corresponding optimization problem. Then, check the  performance of all possible arrangements with the resulting trajectories, and choose the best configuration.

Assuming identical agents, the true anomaly  of the first agent ($\nu$) in orbit $o$(agent $j_1$) should vary between 0 and $\frac{360}{m_o}^{\circ}$. Also, the following constraints should be added to the optimization problem for adjacent agents in orbit $o$ (Figure \ref{fg:DistributeSats}).

\begin{equation}
\label{eq:SatsDisturbute}
\left\lbrace
	\begin{array}{l}
		\nu _{j_2}-\nu _{j_1} =\frac{360}{m_o}^{\circ} \\
		\nu_{j_3}-\nu_{j_2}=\frac{360}{m_o}^{\circ}\\
		\vdots  \\
    	\nu_{j_{m_o}}-\nu_{j_{m_o-1}}=\frac{360}{m_o}^{\circ}\\
	\end{array}
\right.
\end{equation}


\subsection{Collision Avoidance Constraint}

To avoid collision between agents and moving objects around the Earth, the optimization problem should guarantee that all the sensors retain a minimum safe distance from other moving objects during the mission. This can be incorporated in the optimization problem (Equation (\ref{eq:GeneralOptimization})) by heavily penalizing the cost function whenever the minimum distance ($d$) between agents and any orbiting objects in all time instances  is less than an acceptable threshold $tr$ (Equation (\ref{eq:PenalizingCostFunction})).

\begin{equation}
\label{eq:PenalizingCostFunction}
f(p)=\left\{ \begin{array}{l l}
J_{sum}(p) & d \geq tr \\
0 & d < tr \\
\end{array} \right.
\end{equation}

\subsection{Altitude Constraints}

Satellite altitude is an important parameter that relates to many mission requirements.
In the mission design process, many factors such as launch constraint, type of satellite, desired perception resolution, mission life time, groundtrack repeatability and mission expenses  directly affect the altitude of the satellite.

For instance, regardless of the fuel required to retain orbital speeds, depending upon the material technology, a satellite can withstand a limited amount of heating caused by atmospheric resistance. Therefore, for any SSA mission, there is a minimum operational altitude associated with a satellite. This restriction can be included in the trajectories optimization problem by adding the following non-linear constraint.

\begin{equation}
\label{eq:MinimumAltitude}
a_j(1-e_j)>R_E+\underline{c}_j, \; \; \; \; j=1,\cdots,m
\end{equation}

Where $R_E=6371$ Km is the earth radius, and $\underline{c}_j$ denotes the minimum orbital altitude for agent $j$.

Depending upon the altitude, an orbiting satellite receives different kinds of radiation (such as cosmic rays,  Van Allen radiation, solar flares etc.) which can cause damage to the satellite equipment.
To avoid this hazardous exposure, the mission designer can choose different strategies based on the constraints and key requirements.
For a satellite within LEO, one way to stay away from a region of intense radiation is by keeping the altitude of the satellite less than Allen radiation inner belt altitude which is around 1000 Km. Another way is to use expensive radiation-tolerant components in satellite building process to make sure that equipments serve properly throughout the mission.

Therefore, for an space mission,  an upper bound associated with the altitude might be defined.  The upper bound constraint can be expressed by  Equation (\ref{eq:MaximumAltitudeConstraint}).

\begin{equation}
\label{eq:MaximumAltitudeConstraint}
a_j(1-e_j)<R_E+\bar{c}_j, \; \; \; \; j=1,\cdots,m
\end{equation}
Here, $\bar{c}_j$ denotes the agent $j$'s maximum orbital altitude.

\subsection{Specify orbit type}
A satellite  is  classified into various categories according to its orbit altitude, inclination, eccentricity and period etc.. Practically every constraint regarding the shape of the orbit expressed by orbital parameters can be included in the proposed trajectory optimization problem. For instance, to design an orbit for agent $j$ within Low Earth orbit (LEO), its altitude should vary in the 0-2000 km range, i.e. ($0<a_j<R_E+2000$). Similarly,  the motion of agent $j$ moving in a polar orbit  can be planned by adding the equality constraint  $\iota_j=90 {}^{\circ}$ to the problem.
\subsection{Sun-synchronous orbit}
In science missions, one of the most widely used types of orbits is Sun-Synchronous Orbit (SSO). SSO is a near polar and almost-circular geocentric orbit whose nodal precession rate ($\dot{\Omega}$) is equal to the earth's mean rotation rate around the sun. Geometrically, a SSO approximately orients in such a way that the angle between the orbital plane and the vector from Sun to Earth remain the same during the mission. Therefore, the illumination angle of the groundtrack will be  constantly maintained throughout the mission. Beside this interesting characteristic, a SSO has other orbital properties that makes it highly desirable for various applications \cite{ABC_SSO}.
In the subsequent subsections, the SSO orbital parameters constraints imposed by a common scientific mission are briefly discussed, and it is shown how the mission requirements greatly restrict the feasible space of the optimization parameters.

\subsubsection{Precession Rate Constraint}
Due to the out of plane gravitational force caused by Earth's equatorial bulge, an orbit plane gyroscopically precesses. The corresponding nodal precession rate $\dot{\Omega}$  can be operatively computed by Equation (\ref{eq:PrecessionRate}). To have a SSO, $\dot{\Omega}$  should be equal to the Earth's mean orbital rate around the Sun, i.e. $\dot{\Omega}=\frac{365}{365.242199} degree/ day=1.991063802746144 \times 10^{-7} rad/sec$.

\begin{equation}
\label{eq:PrecessionRate}
\dot{\Omega}=-\frac{3}{2} \tilde{J}_2 {\left( \frac{R_E}{a(1-e^2)}\right) }^2 \sqrt{\frac{\mu}{a^3}} cos(\iota)
\end{equation}

Where $ \tilde{J_2}=1.08263 \times 10^{-3}$ denotes the earth's dimensionless zonal harmonic coefficient.

\begin{figure}[h!]
\begin{center}
\includegraphics[width=.43\textwidth]{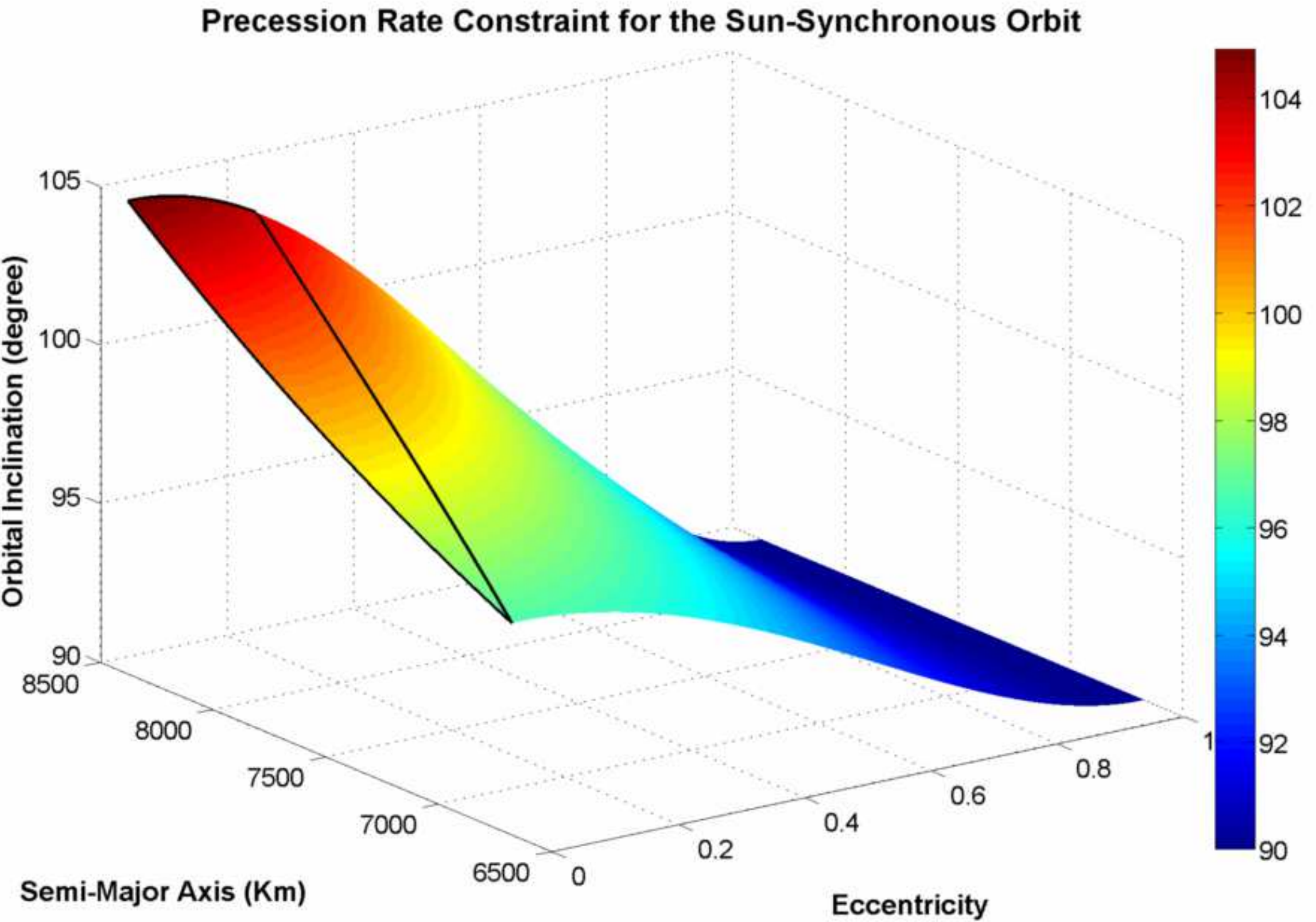}
\caption{Surface of all Sun-synchronous orbits within LEO that satisfy the precession rate condition.}
\label{fg:SunSyncPrecessionRate}
\end{center}
\end{figure}

Figure \ref{fg:SunSyncPrecessionRate} illustrates the surface on which the precession rate condition is satisfied for a SSO within operational region of LEO where $R_E+300Km \leq a \leq R_E +2000 Km$.  Although every point on this surface meets the precession rate condition, relatively a small portion of the surface contains admissible orbits parameters. For instance, most points on the surface don't meet the minimum altitude requirement. In Figure  \ref{fg:SunSyncPrecessionRate}, the area surrounded by black lines contains the points satisfying the minimum altitude constraint $a(1-e)>R_E+300Km$.

\subsubsection{Frozen Orbit Constraint}
Owing to the perturbing forces caused by the oblateness of the earth, the satellite nominal trajectory based on the selected orbital elements does not remain fixed, but it varies as a function of time. In SSO orbit, this fact is exploited to achieve the gyroscopic precession of the orbital plane.
On the other hand, perturbing forces also adversely affect the eccentricity,  $e$, and argument of the perigee, $\omega$, of the orbit.
Perturbation theory states that systematic choices of the orbital parameters can minimize the drift from the selected initial values. In a frozen SSO orbit, parameters $e, \iota$ and $\omega$ are picked in such a way that the secular perturbations of $J_2$ and $J_3$ cancel out, and parameters only undergo a relatively small periodic perturbation with the period equal to the orbit period (Equation (\ref{eq:FrozenOrbit})). Theoretically, a satellite in such an orbit requires a minimum propellant usage during a long term mission.

\begin{equation}
\label{eq:FrozenOrbit}
\left\lbrace \begin{array}{l}
\tilde{\omega} = 90 \; \text{or} \;  270 \;  \;  \; \text{degrees} \\
\tilde{e}=-\frac{J_2}{J_3}\frac{ sin(\iota)}{2  a (1-e^2)}\\
\end{array} \right.
\end{equation}

Where $\tilde{\omega}$ and  $\tilde{e}$ are the mean value of desired frozen orbit parameters.

\begin{figure}[h!]
\begin{center}
\includegraphics[width=.45\textwidth]{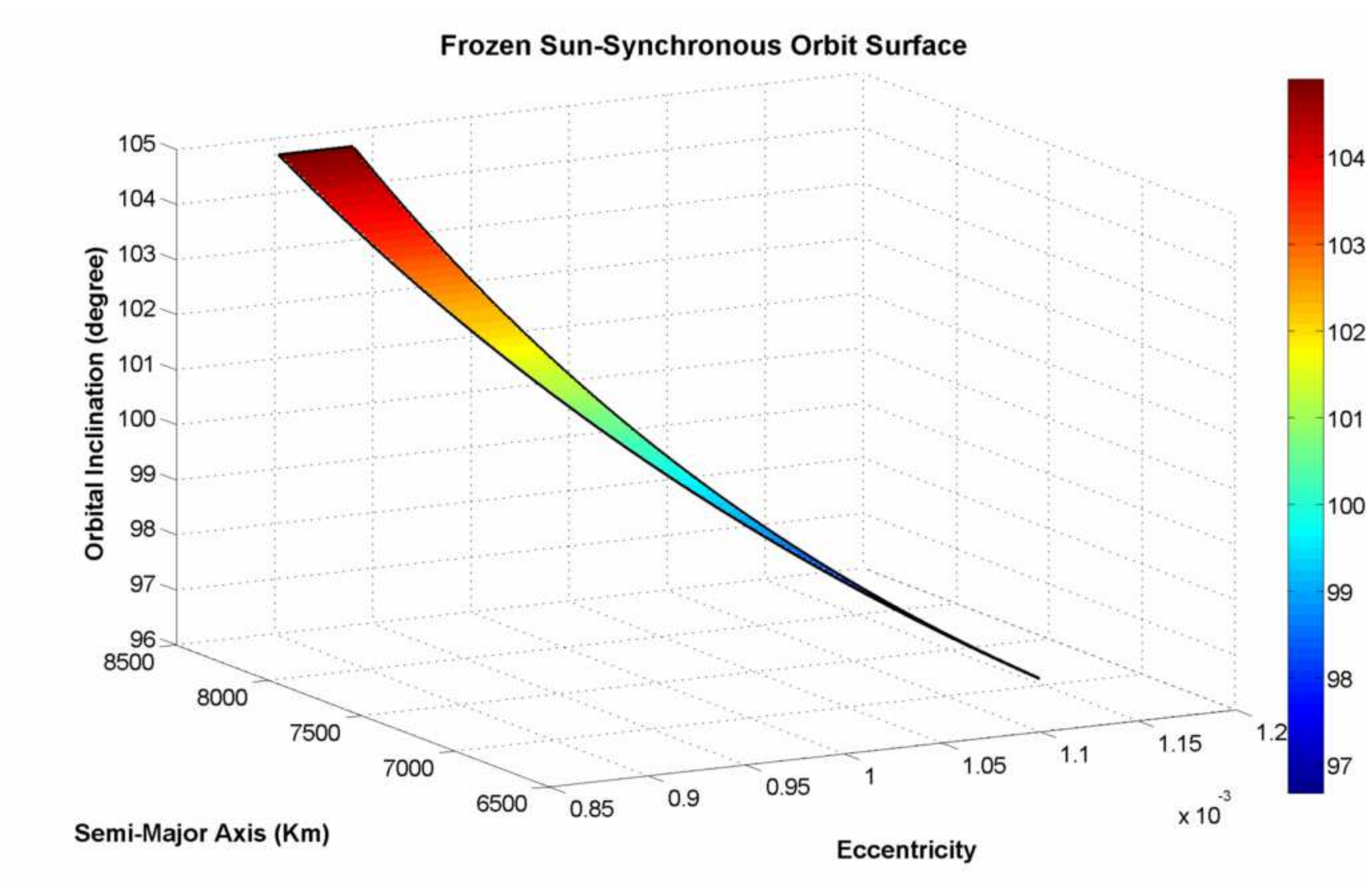}
\caption{The locus of frozen Sun-synchronous orbits parameters within LEO.}
\label{fg:FrozenSunSyncSurface}
\end{center}
\end{figure}

Figure \ref{fg:FrozenSunSyncSurface} depicts the locus of the all Sun-synchronous orbit parameters that satisfy the frozen orbit constraint. As highlighted,  a very small part of the precession rate surface results in the minimum perturbed orbital parameters set.  This has a direct practical implication for the optimization algorithm performance by drastically decreasing the size of the feasible sets.

\section{Numerical Optimization Algorithm}
\label{sec:NumericalOptAlg}
This section considers the problem of choosing an adequate optimization strategy to solve the trajectory design problem governed by Equation (\ref{eq:GeneralOptimization}).
Optimization algorithm selection is the most important part of solving a practical optimization problem.  This process is performed by considering different factors that can affect the overall performance. The performance of an algorithm is commonly assessed according to its convergence speed, solution optimality, robustness to perturbation and uncertainty, required resources, computational complexity and implementation difficulty etc.. Quantifying all the effective factors  is really difficult and occasionally impossible for an optimization method without having knowledge about the problem itself.

Equation (\ref{eq:GeneralOptimization}) is a deterministic constrained non-linear optimization problem with continuous variables. For such a large-scale system with several targets and agents,  the idea of utilizing global optimization methods seems contrived. In fact, current known global approaches offer successful performance in small and medium size problems, as the solution space size exponentially increases with respect to the decision variables number.  Therefore, in this paper, a local optimization algorithms is employed to tackle (\ref{eq:GeneralOptimization}).

In practice, variations of  Sequential Quadratic Programming  (SQP) and  Interior Point Method(IPM) are commonly-used algorithms to solve large-scale general constrained optimization problems. In each iteration, a SQP method generates admissible steps toward local minimum by solving a quadratic model of the objective  function subject to a linearized version of the  problem constraints.
SQP type algorithms have been quite successful  especially in dealing with nonlinear constrains, and has shown more robustness to badly scaled optimization problem than IPM algorithms.

Like SQPs , IPM type algorithms are also well-known for their superior performance in solving nonlinear constrained programs. IPM promising performance in solving linear programming motivated experts to utilize its key ideas, primal-dual steps, to  devise powerful nonlinear optimization algorithms.  IPMs often outperform SQPs in dealing with large-scale applications especially when the system has a block-structure and sparse representation. On the other hand, IPMs have shown lack of robustness to initial point selection, problem scaling and barrier parameters.  Recently, successful attempts have been made to develop more robust IPM-type algorithms \cite{nocedal1999numerical}.

In this paper, since the trajectory design problem formulation is block-structured in terms of optimization parameters, to treat large-scale optimization problems, IPM-type solvers are employed, and for small-scale problems, SQP methods are used.

\section{Simulations}
\label{sec:Sim}
To overcome obstacles in ground-based SSA in detecting and tracking space objects, Spaced-Based SSA (SBSSA) has begun to build a network of satellites equipped with different sensors. To fulfill this purpose, on 26 September 2010, \textit{Minotaur IV} performed its first orbital launch of the Spaced Based Space Surveillance (SBSS) system. The main objective of the SBSS program is to search, detect and track objects in Earth orbits, especially geosynchronous orbit (GEO) objects. SBSS program is about to build a larger constellation of observer satellites to cover wider areas of space.  With the development in technology, it is expected to have more SBSSA missions in near future. The simulation part of this paper is intended to exploit the proposed multi-agent trajectory design method to plan agents movement in the context of SBSSA.

This section is comprised of two parts. In the first subsection, an illustrative simulation study for a fairly small-scale system with five RSOs and two agents is presented. In this system, all targets and agents are orbiting within LEO and in the equatorial plane. In section \ref{subsec:MedLargeScaleSim}, multi-agent orbit design
is performed for more general scenarios of unperturbed SSO.
\begin{figure}[!htb]
\begin{center}
\includegraphics[width=.5\textwidth]{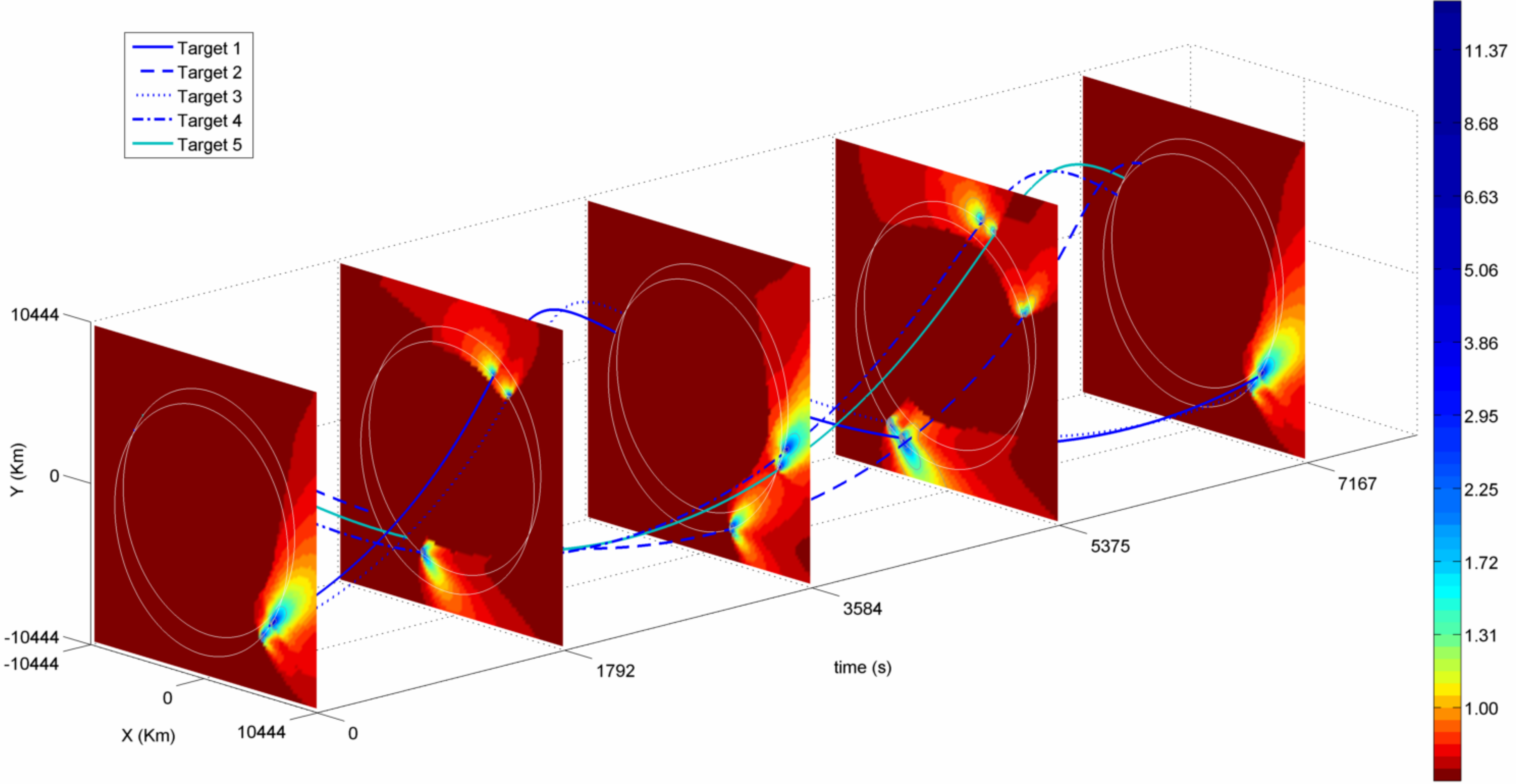}
\caption{ Observation quality field evolution in time.}
\label{fg:SimStudy2DSlices}
\end{center}
\end{figure}
To perform numerical integration of motion dynamic of all objects (Equation (\ref{eq:StateEquation})), the fourth order Runge-Kutta (RK4) method is employed with an adequate step size.
Since the chosen method to solve the trajectory optimization problem results in the largest function value in some feasible neighbourhood, i.e. local optimum solution, to achieve a decent sub-optimal solution for each problem, the optimization algorithm will be run from different initial parameters.


\subsection{Illustrative Case Study}
\label{SimSubsec:2DCase}
In this case study, five RSOs are moving in two elliptic orbits in the equatorial plane with the following Keplerian elements at time $t=0$.

\begin{table}[!htb]

\begin{center}
\begin{tabular}{|l|c|c|c|c|c|c|}
\hline
&$\mathbf{a}$ & $\mathbf{e}$&$\mathbf{\iota} $&$\mathbf{\Omega}$&$\mathbf{\omega}$  & $\mathbf{\nu}$  \\   \hline
{Orbit 1}&8033.72&0.126&0.0&0.0&68.02 &  247.2,  62.2,   237.2 \\  \hline
{Orbit 2}&7898.35&0.057&0.0&0.0&225.2 &  280.8, 270.88\\  \hline
\end{tabular}
\end{center}
\end{table}

The objective is to plan two agents motion so that the visual information gathered from two faces of the targets is maximized. The simulation duration is equal to five times of the  period of the slowest orbit($t_2-t_1=5 \times 2\pi\sqrt{\frac{a^3}{\mu}}=35830.7 \; Sec $), and the sunlight direction is assumed to be from positive sides of $x$-axis to its negative side, i.e. $[-1 \; 0 \; 0]^T$.  The threshold for collision avoidance constraint is considered as $tr=3 Km$.

Figure \ref{fg:SimStudy2DSlices} depicts how the observation quality field changes in five time instances as the system evolves. As shown, the quality field varies in highly non-linear manner.
When the sun shines on one side of the Earth, it casts a shadow on the opposite side of the Earth.

To provide a basis for comparison, two different scenarios are considered. In the first situation, it is assumed that the agents are co-orbital and uniformly distributed on an equatorial orbit. While in the second case, two different equatorial orbits are designed. For both cases, twelve unknown variables should be determined by solving the optimization problem, i.e. $ p=[a_1,e_1,\Omega_1,\omega_1,\iota_1,\nu_1,a_2,e_2,\Omega_2,\omega_2,\iota_2,\nu_2]^T$. The $p$ vector will be limited to upper and lower  $\bar{p}$ and $\underline{p}$ vectors in Equation $\ref{eq:ScalingVectors2D}$. Since the systems are small-scale in these scenarios, a SQP-type algorithm is employed to numerically solve the corresponding optimization problem.

\begin{equation}
\label{eq:ScalingVectors2D}
\begin{array}{l}
\bar{p}=[8033,0.15,360,360,10,360,8033,0.15,\\ \quad \quad \quad \quad \quad \quad \quad \quad  \quad \quad \quad \quad  360,360,10,360]^T\\
\underline{p}=[7898.3,0,0,0,0,0,7898.3,0,0,0,0,0]^T
\end{array}
\end{equation}

 \subsubsection{Co-orbital Case}
Equation (\ref{eq:Sim2DScenario1}) is the  non-linear optimization problem corresponding to the first scenario.
Assuming the minimum operational altitude for both agents to be 300 Km, the altitude constraints are $a_j(1-e_j)>R_E+300$ for $j=1,2$.  Equality constraints $\Omega_1=\Omega_2=\iota_1=\iota_2=0$ force the orbits to be in the equatorial plane. Constraint $A^T p= 0 $ makes the two agents co-orbital, and constraints  $\nu_1\leq 180$ and $\nu_2-\nu_1= 180$ result in equally spaced agents in the orbit.

\begin{equation}
\label{eq:Sim2DScenario1}
  \left\{
  \begin{array}{l l}
   \underset{p}{ \text{maximize}} &  min( f(p))\\
    \text{subject to :} &  \\
& \underline{p}\leq p \leq \bar{p}\\
&  a_j(1-e_j)>R_E+300, \quad j=1,2 \\
& \begin{array}{l}
\Omega_1=0 \; \; \rightarrow
[ 0_{1 \times 2} \quad 1 \quad 0_{1 \times 11}
] p=0 \\
\Omega_2=0 \; \; \rightarrow
[0_{1 \times 8} \quad 1 \quad 0_{1 \times 3}] p=0
\end{array} \\
&  \begin{array}{l}
\iota_1=0 \; \; \rightarrow
[ 0_{1 \times 4} \quad 1 \quad 0_{1 \times 7}
] p=0 \\
\iota_2=0 \; \; \rightarrow
[0_{1 \times 10} \quad 1 \quad 0] p=0
\end{array} \\
&  A^T p=[ I_{5\times 5} \; 0_{5\times 1} \; -I_{5 \times 5} \; 0_{5\times 1} ] p =0_{5\times 1} \\
& \begin{array}{l}
\nu_1 \leq 180 \rightarrow [ 0_{5 \times1} \quad 1 \quad 0_{5 \times1} \quad 0 ] p \leq 180 \\
\nu_2 -\nu_1 = 180 \rightarrow 
 [ 0_{5 \times1} \; -1 \; 0_{5 \times1} \; 1 ] p=180
\end{array}
\end{array} \right.
 \end{equation}

 \subsubsection{Agents with different orbits}
 The second scenario formulation is like Equation (\ref{eq:Sim2DScenario1}) with only the six first constraints to make sure that the resultant orbits satisfy boundary conditions, minimum altitude constraints, and position in equatorial plane.

\begin{figure}[!h]
\begin{center}
\includegraphics[width=.45\textwidth]{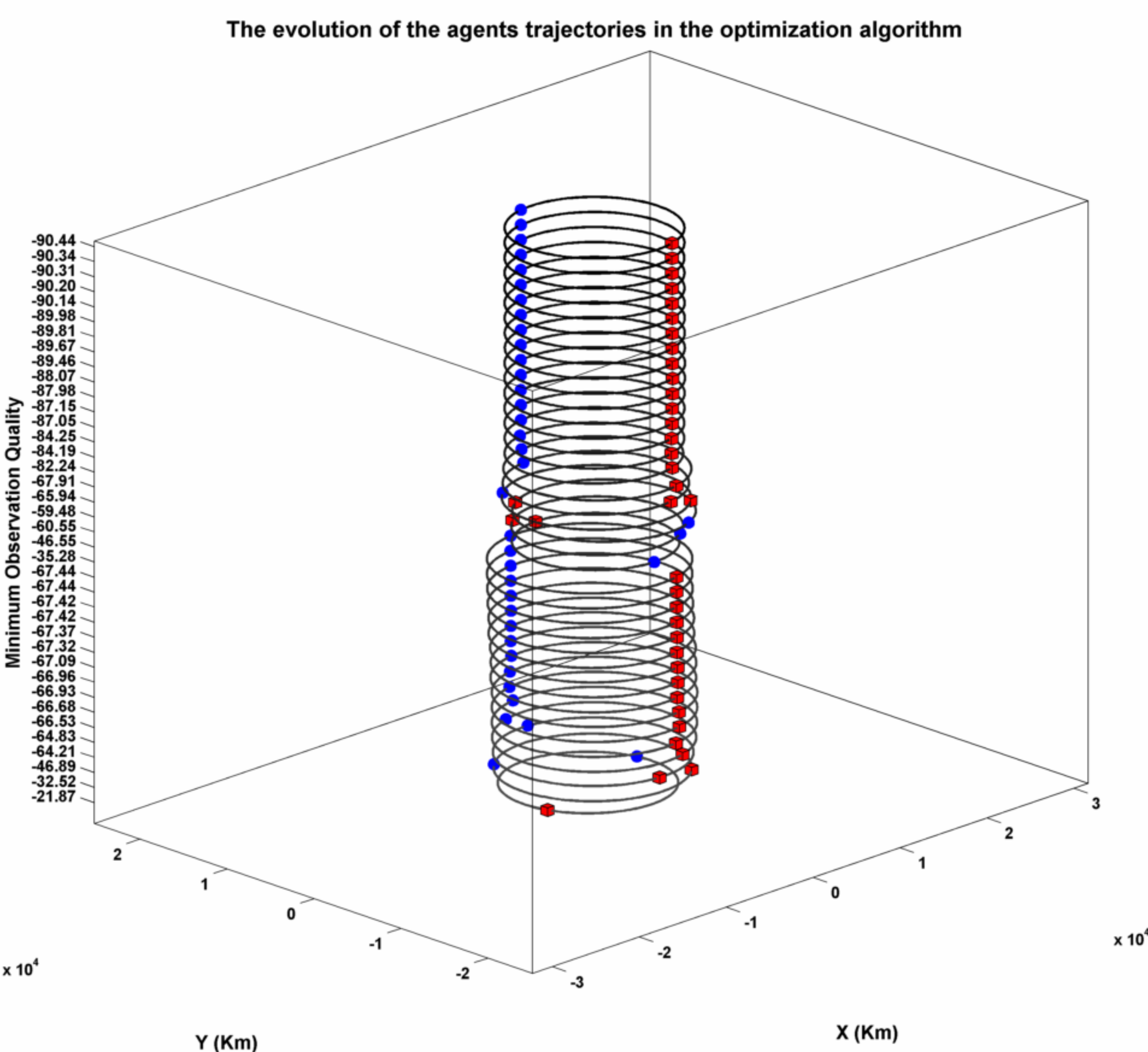}
\caption{The best obtained suboptimal solution of the optimization problem in Equation (\ref{eq:Sim2DScenario1}).}
\label{fg:TrajectoryEvolution1}
\end{center}
\end{figure}

Both optimization problems are run from several feasible initial parameters to obtain satisfactory sub-optimal solutions. Figure \ref{fg:TrajectoryEvolution1} describes how the best achieved trajectory evolves from a feasible initial solution in the co-orbital case. As it appears, the agents are uniformly distributed along the trajectories in all iterations.


\begin{figure}[!htb]
\centering
\includegraphics[width=.45\textwidth]{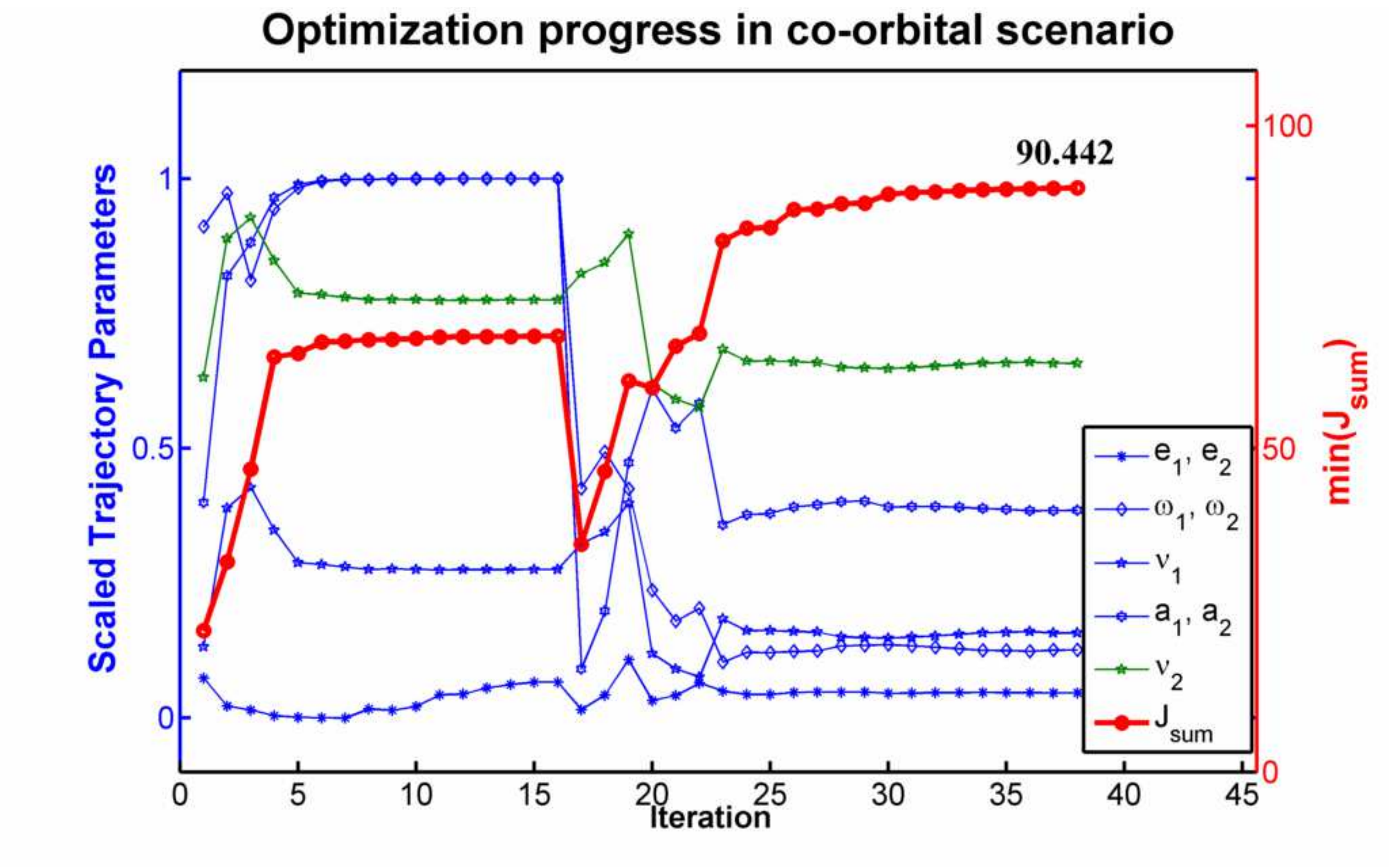} 
\caption{Optimization progress in the co-orbital scenario. }
\label{fg:2DOptOrbitsJSUM_CoOrbital}
\end{figure}

\begin{figure}[!htb]
\centering
\includegraphics[width=.45\textwidth]{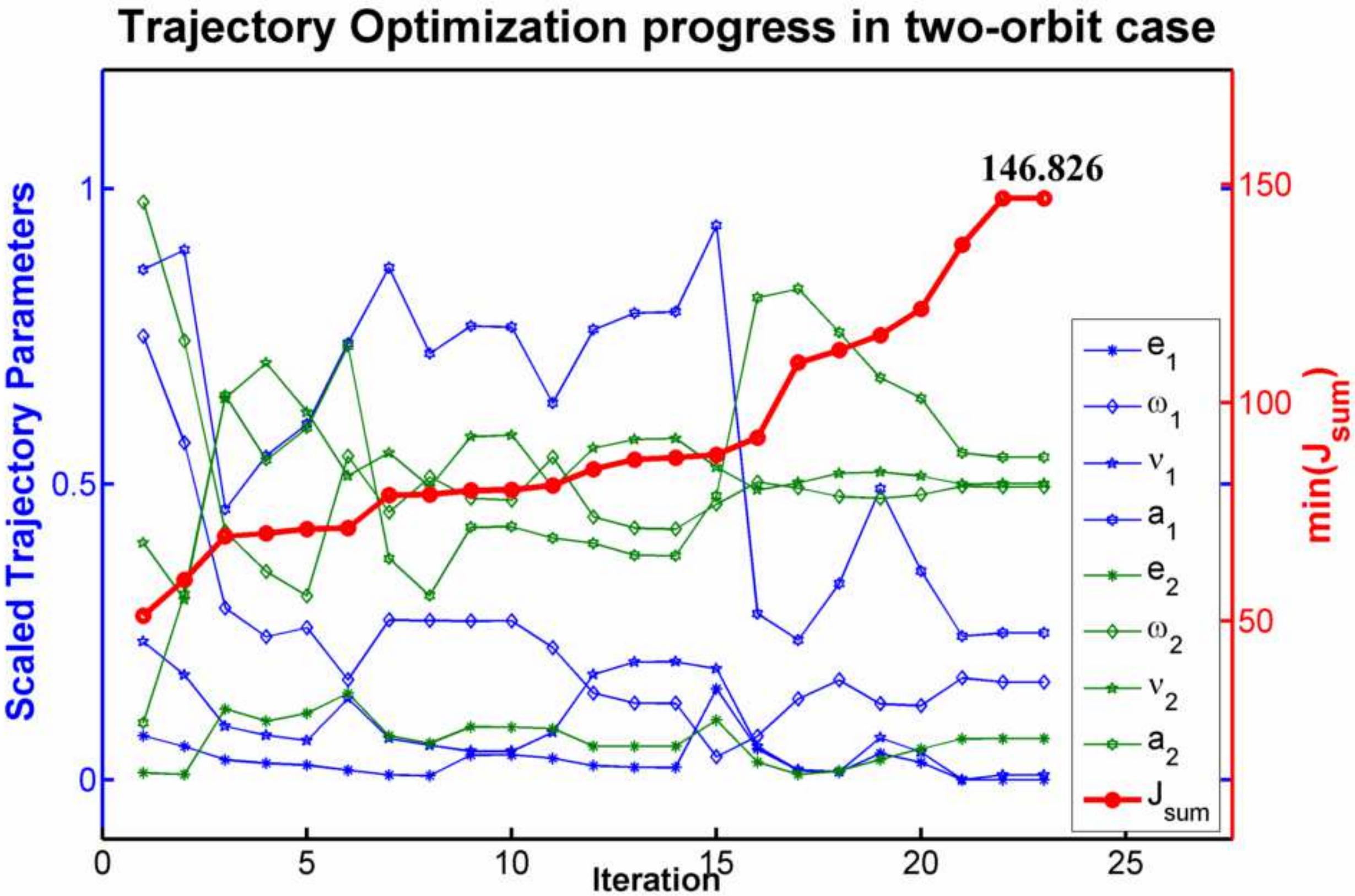} 
\caption{Optimization progress in the scenario with to distinct orbits. }
\label{fg:2DOptOrbitsJSUM_TwoOrbits}
\end{figure}

\begin{figure*}[!htb]
\centering
\begin{tabular}{cc}
\includegraphics[width=0.345\textwidth]{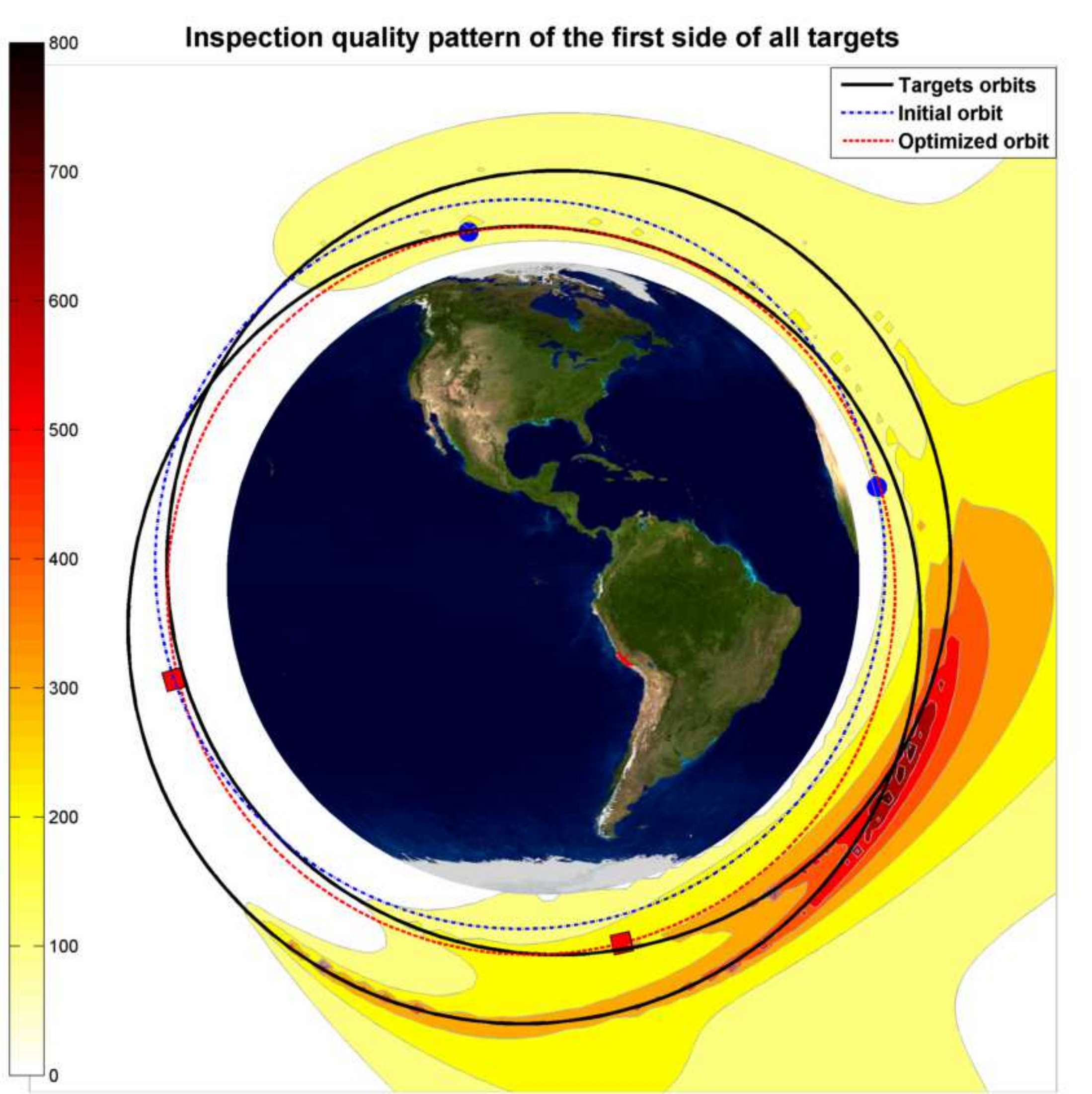} &
\includegraphics[width=0.345\textwidth]{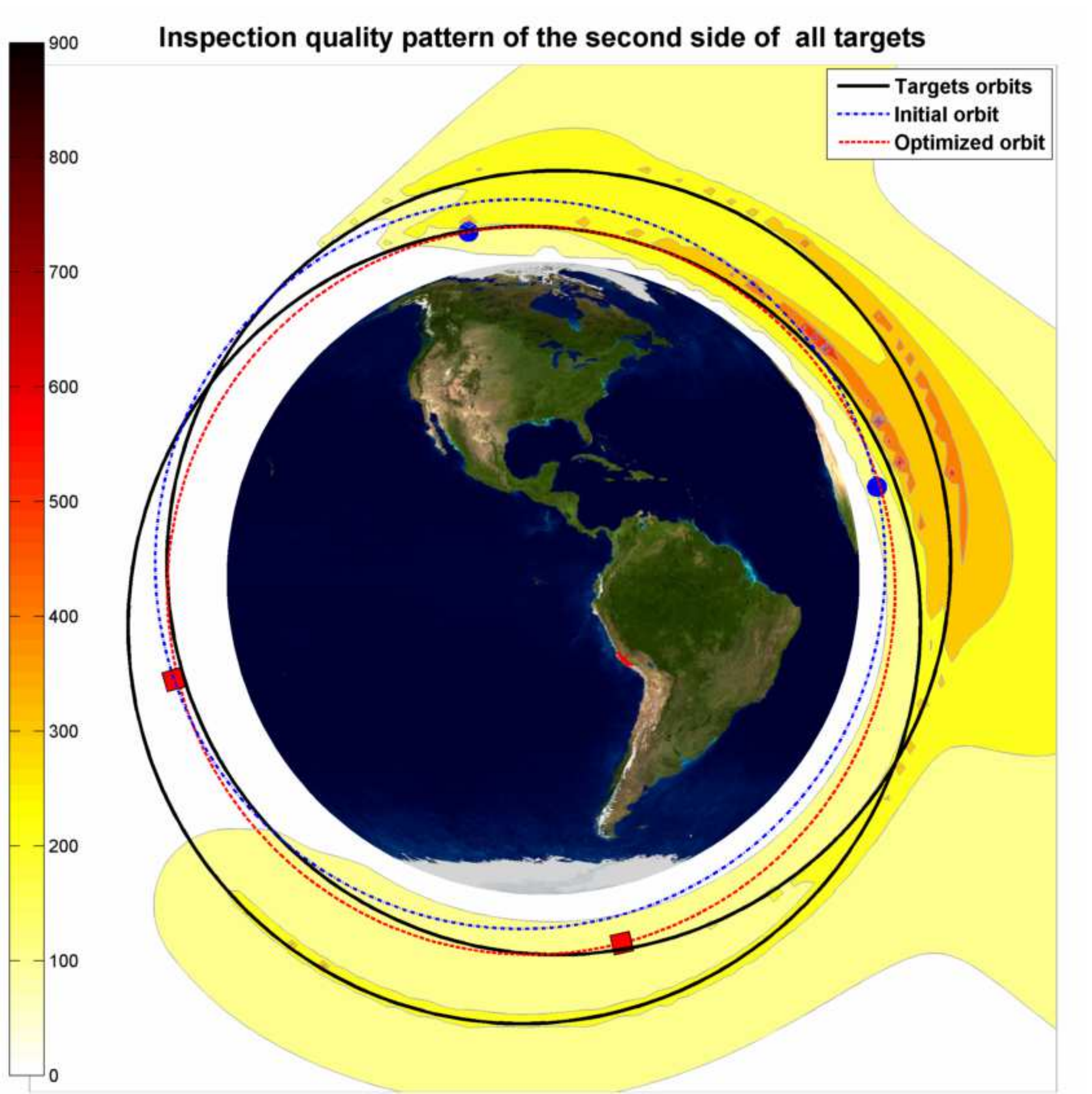} \\
\includegraphics[width=0.345\textwidth]{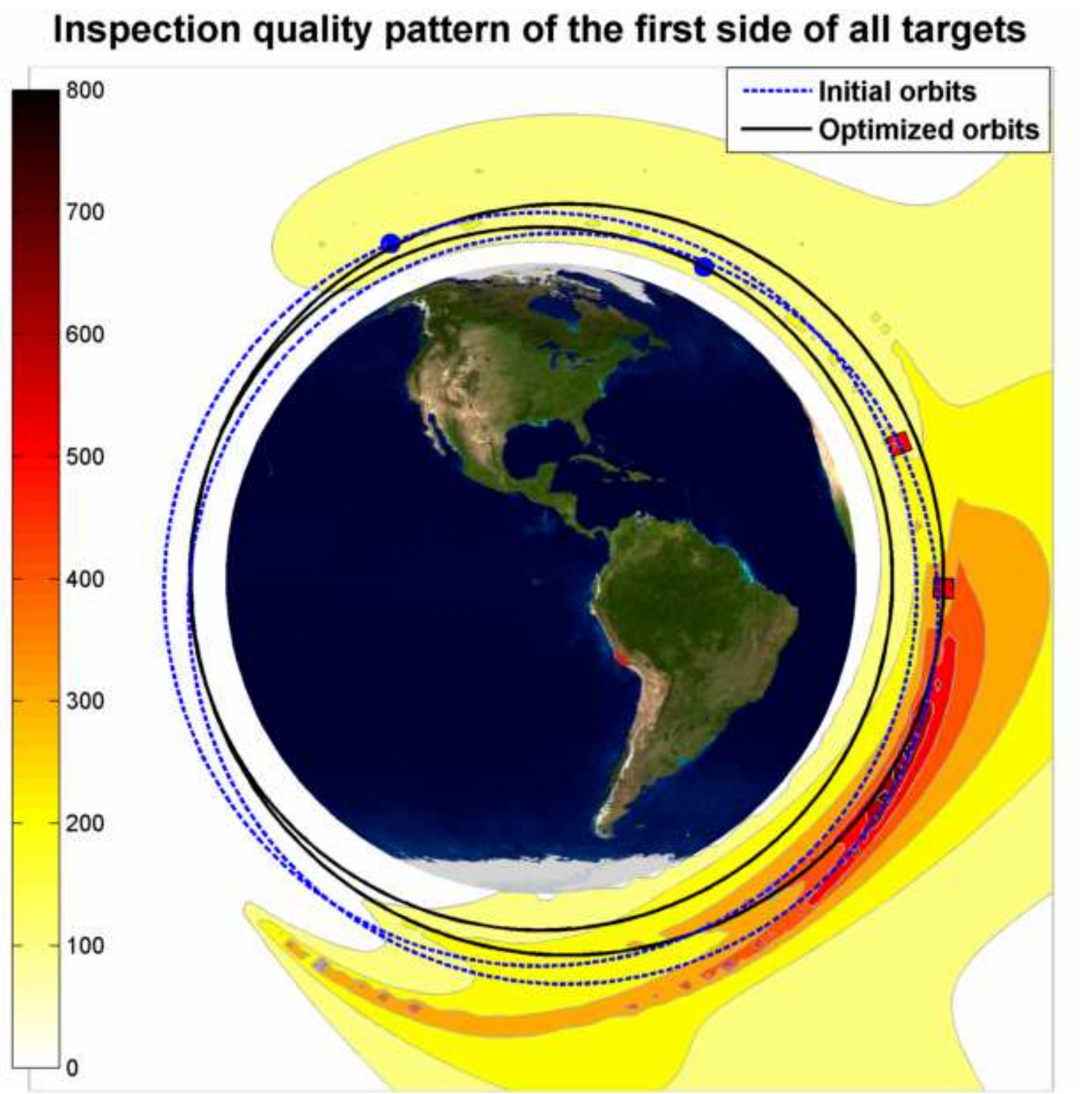} &
\includegraphics[width=0.345\textwidth]{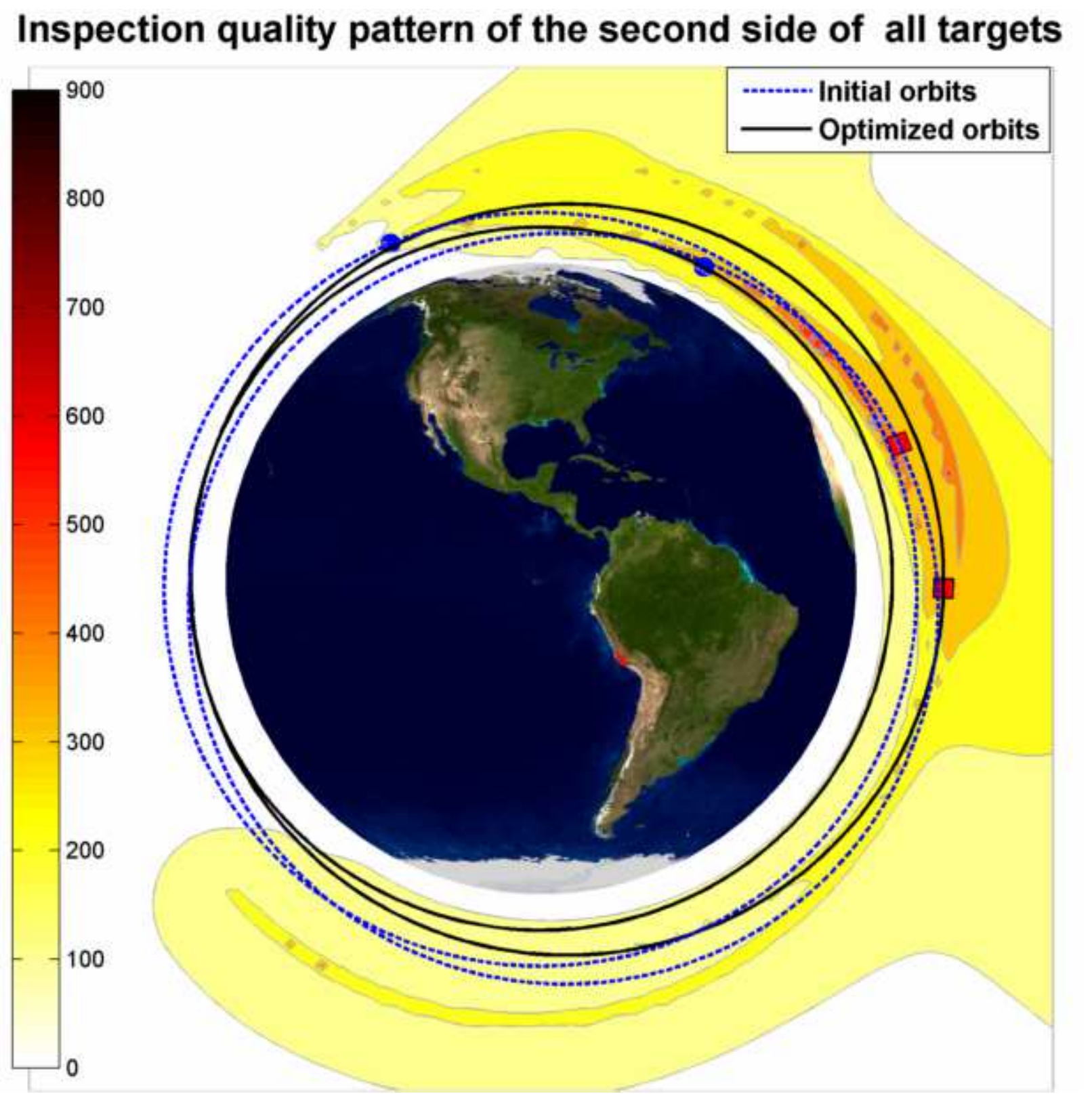} \\
\end{tabular}
\caption{Optimal trajectory and sides quality pattern of all RSOs for both design scenarios. ( Co-orbital case (Top) and two distinct agents trajectories design case (Bottom)).}
\label{fg:2DOptOrbitBothsidesCoOrbital}
\end{figure*}

\begin{figure}[!h]
\centering
\includegraphics[width=0.45\textwidth]{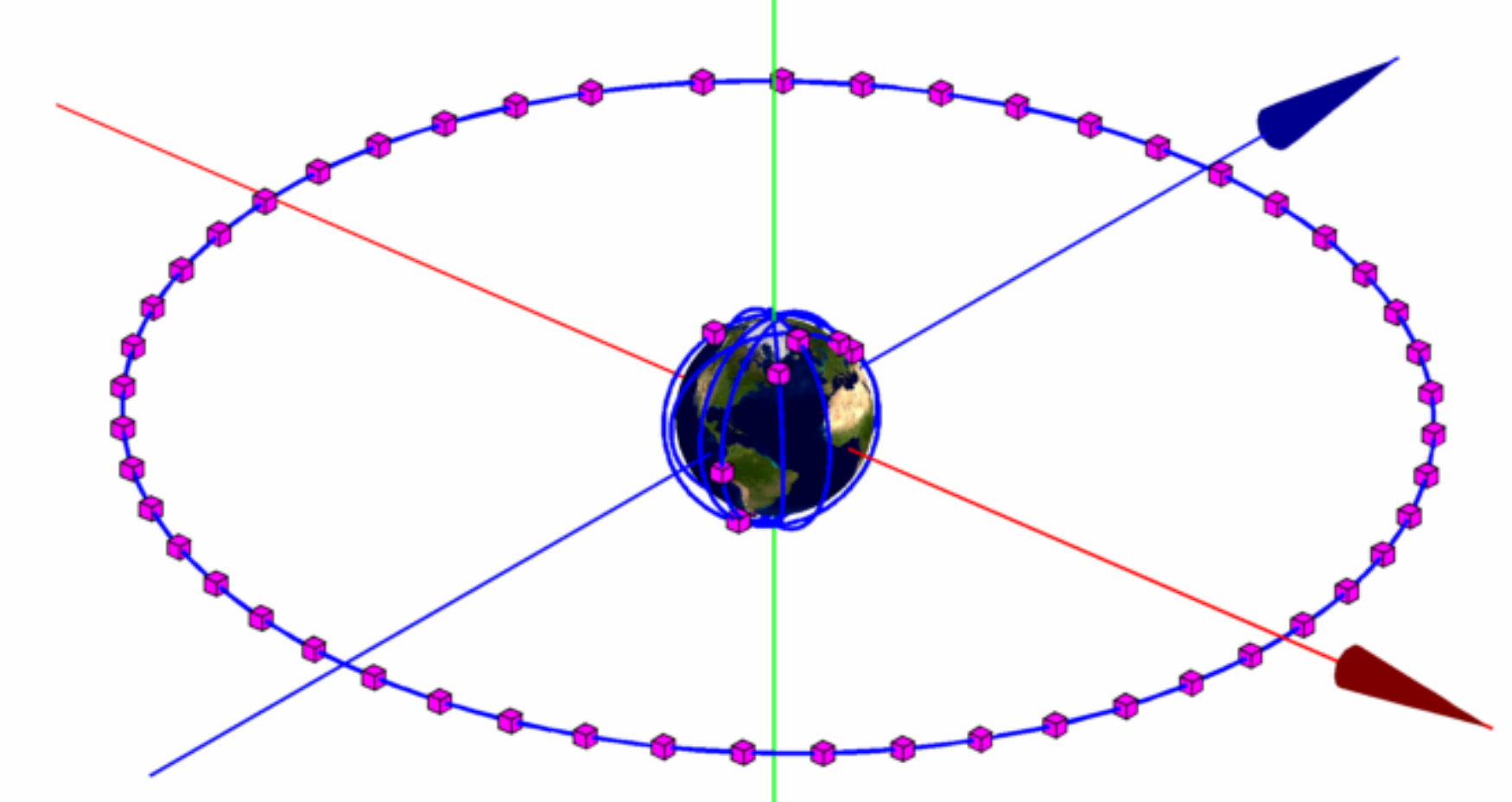}
 \caption{RSOs arrangement in a typical SBSSA application.}
 \label{fg:3DSBSSA_RSOConfig}
\end{figure}

Figure \ref{fg:2DOptOrbitsJSUM_CoOrbital} and \ref{fg:2DOptOrbitsJSUM_TwoOrbits} illustrate how the scaled optimization variables and the minimum observation quality of all targets sides converge for each scenario. As highlighted, after convergence, the optimization cost function in the co-orbital case is 90.449 that is clearly less than the corresponding obtained observation quality along two other trajectories (146.826). This indicates that in systems in which agents have more degree of freedoms, a higher perception level can be achieved.

Figure \ref{fg:2DOptOrbitBothsidesCoOrbital} depicts the sum of the quality pattern for each side of the targets, as well as initial and optimized trajectories for both design scenarios. The graph shows how the initial trajectories are altered and developed to be well positioned in the areas with higher perception level in both sides.

\subsection{Frozen SSO Design}
\label{subsec:MedLargeScaleSim}

In this part, the multi-agent orbit design procedure is performed for larger SBSSA systems. In the first part, because of desirable characteristics of SSO reported for the SBSS program \cite{Yuzhuo2011},  Sun-Synchronous Frozen orbits are found for a SBSSA system with five agents which must characterize sixty RSOs.
In this scenario, fifty one RSOs are traveling on an orbit in GEO and the rest are moving within LEO (Figure \ref{fg:3DSBSSA_RSOConfig}).

In this design problem, the Keplerian parameters of the desired orbits should satisfy SSO and Frozen Orbit constraints (Equation (\ref{eq:PrecessionRate}) and (\ref{eq:FrozenOrbit})). Furthermore, the lower and upper band for altitude is assumed to be $\underline{c}_j =300 Km$ and $\bar{c}_j=1000 Km$,  the orbital inclination lies in $[96.5, 102.5]$ interval and the argument of perigee $\omega_j=90$. Thus the optimization problem can be written as Equation (\ref{eq:Sim3DScenario1}).

\begin{equation}
\label{eq:Sim3DScenario1}
  \left\{
  \begin{array}{l l}
   \underset{p}{ \text{maximize}} &  min( f(p))\\
    \text{subject to :} &  \\
& \underline{p}\leq p \leq \bar{p}\\
&  a_j(1-e_j)>R_E+300, \quad j=1,\cdots,5 \\
&  a_j(1-e_j)<R_E+1000, \quad j=1,\cdots,5 \\
&  -\frac{3}{2} \tilde{J}_2 {\left( \frac{R_E}{a_j(1-e_j^2)}\right) }^2 \sqrt{\frac{\mu}{a_j^3}} cos(\iota_j)=1.99\times 10^{-7}\\
& \quad \quad \quad \quad \quad \quad \quad \quad j=1,\cdots,5\\
& e_j=-\frac{J_2}{J_3}\frac{ sin(\iota_j)}{2  a_j (1-e_j^2)} , \quad j=1,\cdots,5\\
& \omega_j=90  , \quad j=1,\cdots,5
\end{array} \right.
 \end{equation}

\begin{figure}[!htp]
\centering
\includegraphics[width=0.45\textwidth]{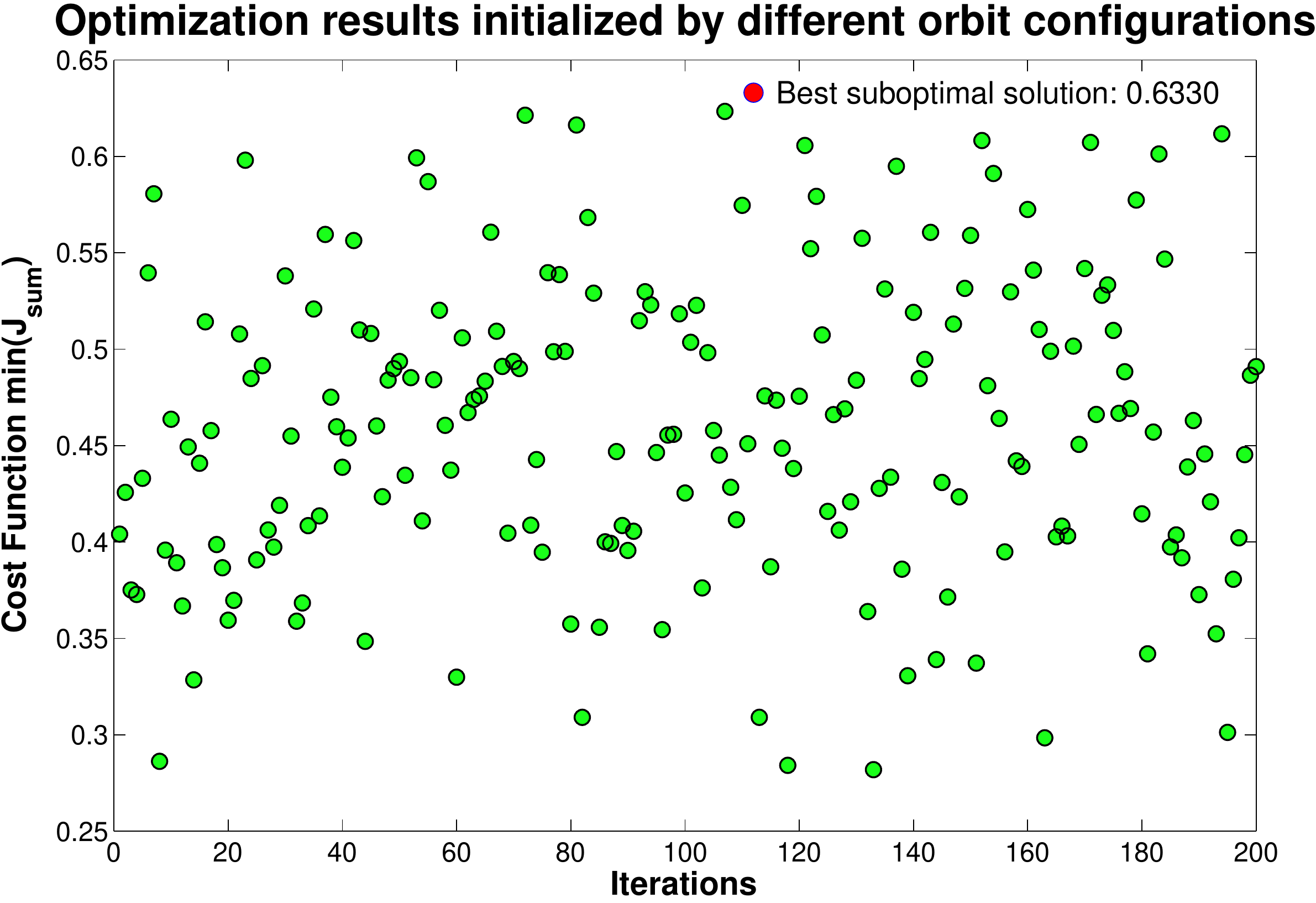}
\caption{The best obtained sub-optimal solution among two-hundred optimization runs for SSO design case.}\label{fg:SunSyncOptGraphs_SubOptimalSol}
\end{figure}

\begin{figure}[!htp]
\centering
\includegraphics[width=0.42\textwidth]{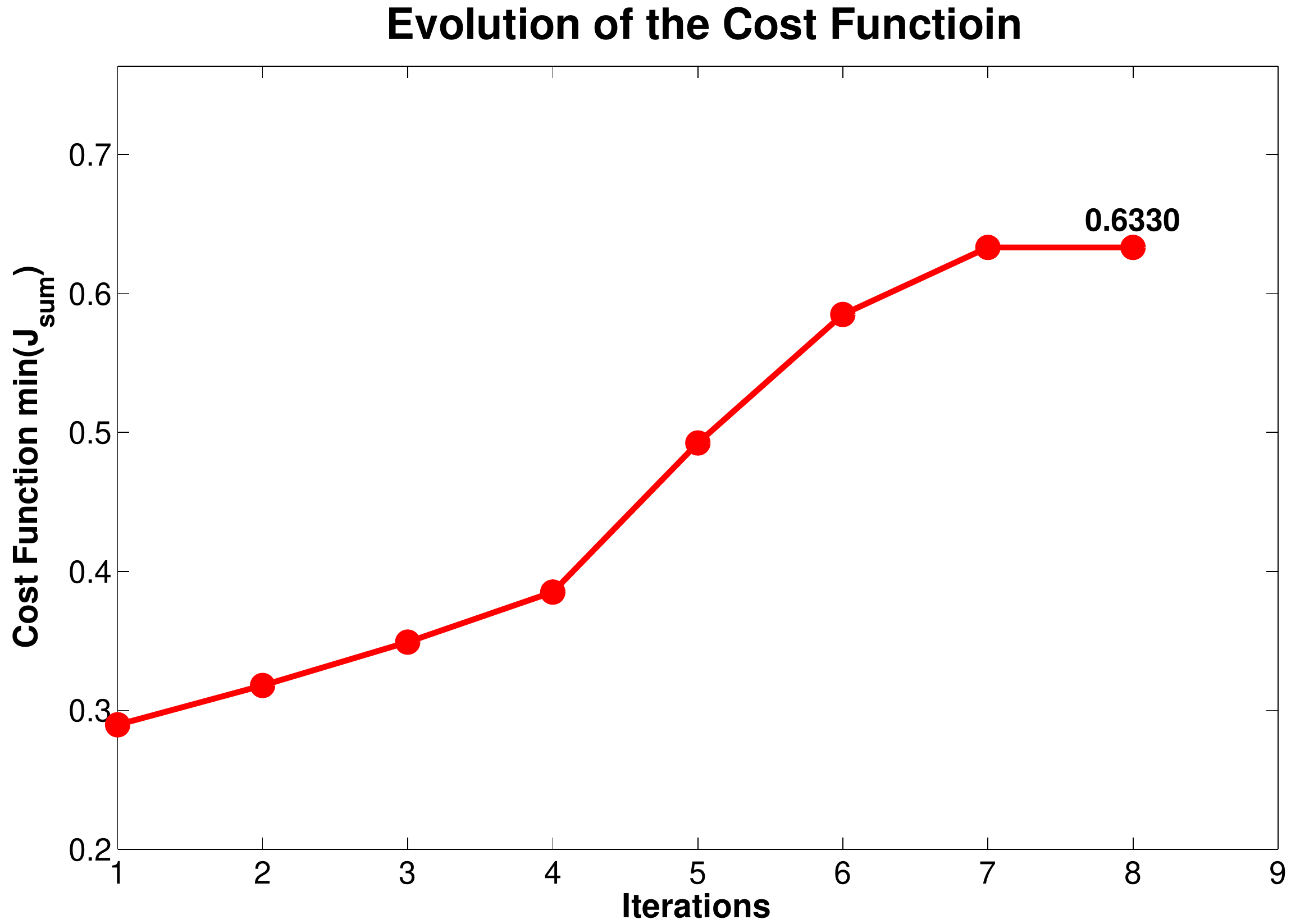} 
\caption{The cost function evolution of the best obtained solution in SSO design scenario. }\label{fg:SunSyncOptGraphs_CostFunc}
\end{figure}

Figure \ref{fg:SunSyncOptGraphs_SubOptimalSol} and \ref{fg:SunSyncOptGraphs_CostFunc} report the  procedure to achieve a  sup-optimal solution for the optimization problem in Equation (\ref{eq:Sim3DScenario1}).  Figure \ref{fg:SunSyncOptGraphs_SubOptimalSol} depicts 200  runs of the optimization from different initial parameters. The best minimum inspection quality for the desired sides of RSOs is 0.633, and Figure \ref{fg:SunSyncOptGraphs_CostFunc} indicates how the initial objective function is evolved in each optimization iteration.

\begin{figure}[!htb]
\centering
\includegraphics[width=0.4\textwidth]{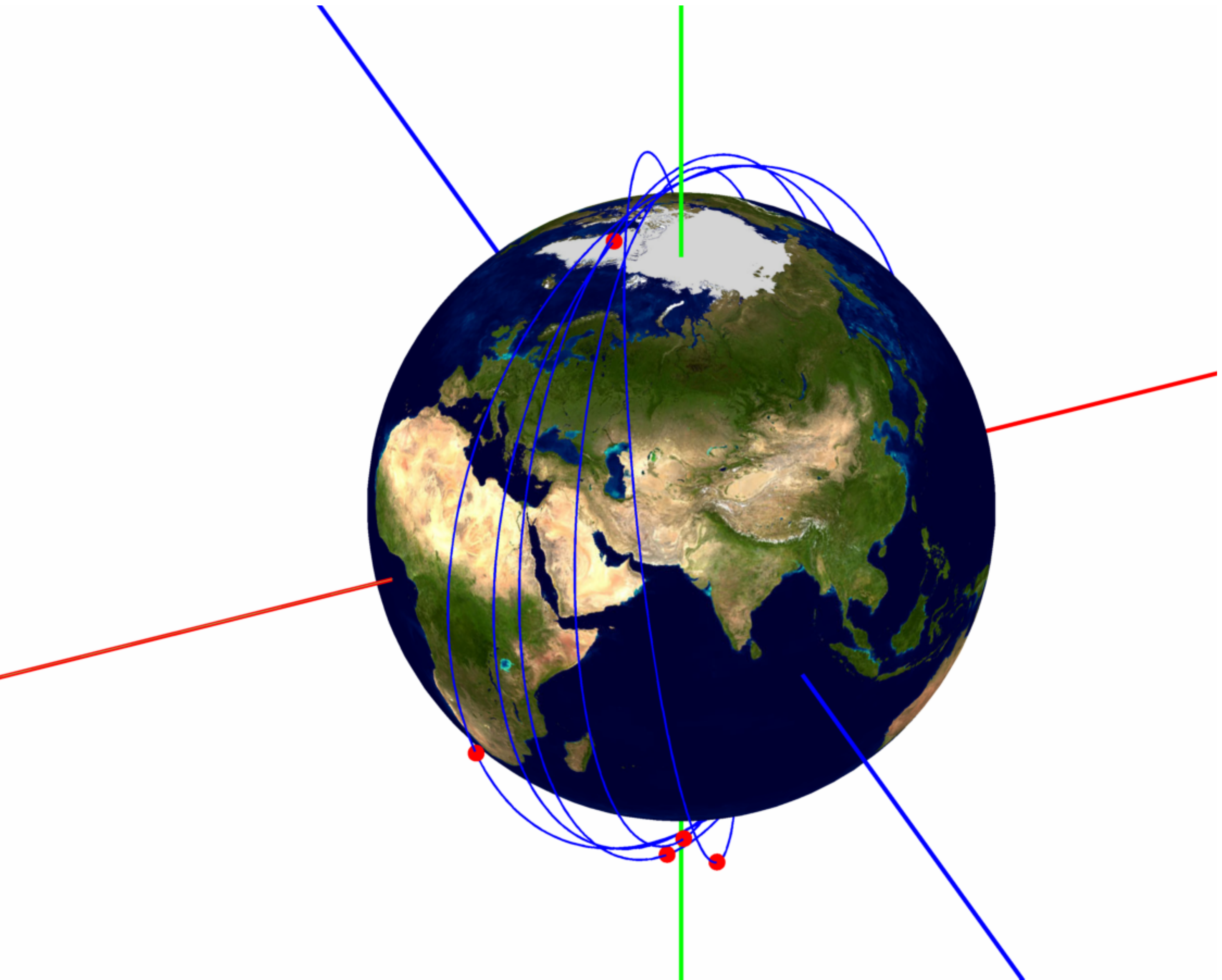}
 \caption{The optimized frozen SSOs for a multi-agent system.}
 \label{fg:3DSSO_FrozenObsPlot}
\end{figure}

Figure \ref{fg:3DSSO_FrozenObsPlot} illustrates the designed trajectories on which the minimum perception quality is maximized, and it corresponds to the best acquired sub-optimal solution.

\section{Conclusion}
\label{sec:Conslution}

In this paper, a systematic way to optimize trajectories for a sensor network system for perception enhancement purpose is presented. This has a major impact on the performance of the resource allocation stage in any sensor management system. Specifically, the proposed method is of interest for large-scale systems characterized by highly dynamic behaviour because of it robustness properties to  possible perturbations and model uncertainties.
Case study results for Spaced-Based Space Situational Awareness (SBSSA) applications indicate the effectiveness of the method in improving received perception qualities along the designed trajectories where the agents cooperatively meet the mission requirements, and less effort is needed to serve the main purpose.

\ifCLASSOPTIONcaptionsoff
  \newpage
\fi



%
%
%

\bibliographystyle{ieeetr}
\bibliography{OrbitDesign_REF}

\begin{thebibliography}{10}

\bibitem{SensorNetVladimirova2}
T.~Vladimirova, C.~Bridges, G.~Prassinos, X.~Wu, K.~Sidibeh, D.~Barnhart, A.-H.
  Jallad, J.~Paul, V.~Lappas, A.~Baker, K.~Maynard, and R.~Magness,
  ``Characterising wireless sensor motes for space applications,'' in {\em
  Adaptive Hardware and Systems, 2007. AHS 2007. Second NASA/ESA Conference
  on}, pp.~43 --50, aug. 2007.

\bibitem{ABC_SSO}
R.~J. Boain, ``A-b-cs of sun-synchronous orbit mission design,'' in {\em
  AAS/AIAA Space Flight Mechanics Meeting}, (Maui, HI; United States), 2004.

\bibitem{orbitDesignGEOSYNC}
E.~Hernandez, J.~Bolivar, and Q.~Wang, ``Geosynchronous transfer orbit design:
  A practical approach,'' in {\em Modeling, Simulation and Applied Optimization
  (ICMSAO), 2011 4th International Conference on}, pp.~1 --4, april 2011.

\bibitem{OrbitDesignNonKeplerian}
Z.~Zhu, J.~Yuan, and L.~Zheng, ``Shape-based method for non-keplerian orbit
  design,'' in {\em 12th Biennial International Conference on Engineering,
  Construction, and Operations in Challenging Environments; and Fourth
  NASA/ARO/ASCE Workshop on Granular Materials in Lunar and Martian
  Exploration: Earth and Space 2010}, (Honolulu, Hawaii, United States),
  pp.~1934--1940, March 2010.

\bibitem{ConvexOrbitDesignMit}
M.~J. Tillerson, {\em Orbit Design Convex Optimization Technique}.
\newblock PhD thesis, Massachusetts Institute of Technology.

\bibitem{Inalhan02relativedynamics}
G.~Inalhan, M.~Tillerson, and J.~P. How, ``Relative dynamics and control of
  spacecraft formations in eccentric orbits,'' {\em AIAA Journal of Guidance,
  Control, and Dynamics (0731-5090}, vol.~25, pp.~48--59, 2002.

\bibitem{FormationFlying}
M.~Tillerson and J.~How, ``Formation flying control in eccentric orbits,'' in
  {\em Proceedings of the AIAA Guidance, Navigation, and Control Conference,
  Montreal}.

\bibitem{Hamid2011}
H.~Nourzadeh and J.~E. McInroy, ``Planning the visual measurement of $n$ moving
  objects by $m$ moving cameras, given their relative trajectories,'' in {\em
  IEEE Multi-Conference on Systems and Control}, September 2011.

\bibitem{HamidRobust2013}
H.~Nourzadeh and J.~McInroy, ``Robust visual measurement planning in
  multi-robot systems,'' in {\em Automation Science and Engineering (CASE),
  2013 IEEE International Conference on}, pp.~176--182, 2013.

\bibitem{HamidPatrolling2013}
H.~Nourzadeh and J.~McInroy, ``Integrated planning of constraint sensor
  management and patrolling,'' in {\em Automation Science and Engineering
  (CASE), 2013 IEEE International Conference on}, pp.~837--843, 2013.

\bibitem{irv97a}
J.~Irvine, ``National imagery interpretability rating scales (niirs): Overview
  and methodology,'' in {\em Proceedings of the International Society for
  Optical Engineering (SPIE)}, vol.~3128, pp.~93--103, July 1997.

\bibitem{irv03a}
J.~Irvine, ``National imagery intelligence rating scale {(NIIRS)},'' in {\em
  The Encyclopedia of Optical Engineering} (R.~Driggers, ed.), Marcel Dekker,
  2003.

\bibitem{lea96a}
J.~Leachtenauer, ``National imagery interpretability rating scales: Overview
  and product description,'' in {\em Proceedings of the American Society of
  Photogrammetry and Remote Sensing Annual Meetings}, April 1996.

\bibitem{lea97a}
J.~C. Leachtenauer, W.~Malila, J.~M. Irvine, L.~Colburn, and N.~Salvaggio,
  ``General image-quality equation: {GIQE},'' {\em Applied Optics}, vol.~36,
  pp.~8322--8328, 1997.

\bibitem{mav95a}
L.~Maver, C.~Erdman, and K.~Riehl, ``Imagery interpretability rating scales,''
  in {\em Digest of Technical Papers, International Symposium Society for
  Information Display}, vol.~XXVI, pp.~117--120, 1995.

\bibitem{young_signal_2008}
S.~S. Young, R.~G. Driggers, and E.~L. Jacobs, {\em Signal Processing and
  Performance Analysis for Imaging Systems}, pp.~16--25.
\newblock Artech House Publishers, Apr. 2008.

\bibitem{ModulationTransferFunction}
G.~D. Boreman, {\em Modulation Transfer Function in Optical and
  {ElectroOptical} Systems}, pp.~31--42.
\newblock {SPIE} Publications, July 2001.

\bibitem{hecht2002optics}
E.~Hecht, {\em Optics}, pp.~347--472.
\newblock Pearson education, Addison-Wesley, 2002.

\bibitem{forsyth_computer_2002}
D.~A. Forsyth and J.~Ponce, {\em Computer Vision: A Modern Approach},
  pp.~9--11.
\newblock Prentice Hall, {US} ed~ed., Aug. 2002.

\bibitem{chen_active_2008}
S.~Chen, Y.~F. Li, J.~Zhang, and W.~Wang, {\em Active Sensor Planning for
  Multiview Vision Tasks}, pp.~4--5.
\newblock Springer, 1~ed., 2008.

\bibitem{curtis_orbital_2005}
H.~Curtis, {\em Orbital Mechanics: For Engineering Students}, pp.~158--161.
\newblock {Butterworth-Heinemann}, Jan. 2005.

\bibitem{nocedal1999numerical}
J.~Nocedal and S.~Wright, {\em Numerical Optimization}, pp.~352--417.
\newblock Springer Series in Operations Research, Springer, 1999.

\bibitem{Yuzhuo2011}
S.~Z. Hou~Yuzhuo, Huang~Xuexiang, ``A satellite orbit design method for
  space-based space surveillance,'' in {\em European Space Surveillance
  Conference}, June 7-9 2011.

\end{thebibliography}

%

\begin{IEEEbiography}[{\includegraphics[width=1in,height=1.25in,clip,keepaspectratio]{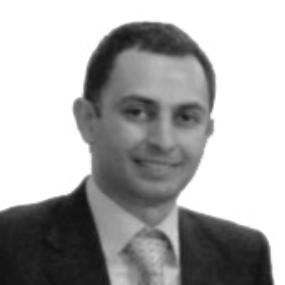}}]{Hamidreza  Nourzadeh}

was born in Tehran, Iran, in 1978. He received the B.S.
degree from Shahid Beheshti University, Tehran, Iran, in
2003, the M.S. degree from School of Electrical
Engineering, K. N. Toosi University of technology,
Tehran, Iran, in 2006, and the Ph.D. degree from University of Wyoming, Laramie, Wyoming, in 2013, all in Electrical Engineering.
He is currently a postdoctoral research associate at Electrical Computer and System Engineering Department 
at the Rensselaer Polytechnic Institute, NY,
USA.

His interests are in the areas of control systems, multi-agent systems, 3-D Computer Vision, optimization and  system identification.
\end{IEEEbiography}

\begin{IEEEbiography}[{\includegraphics[width=1in,height=1.25in,clip,keepaspectratio]{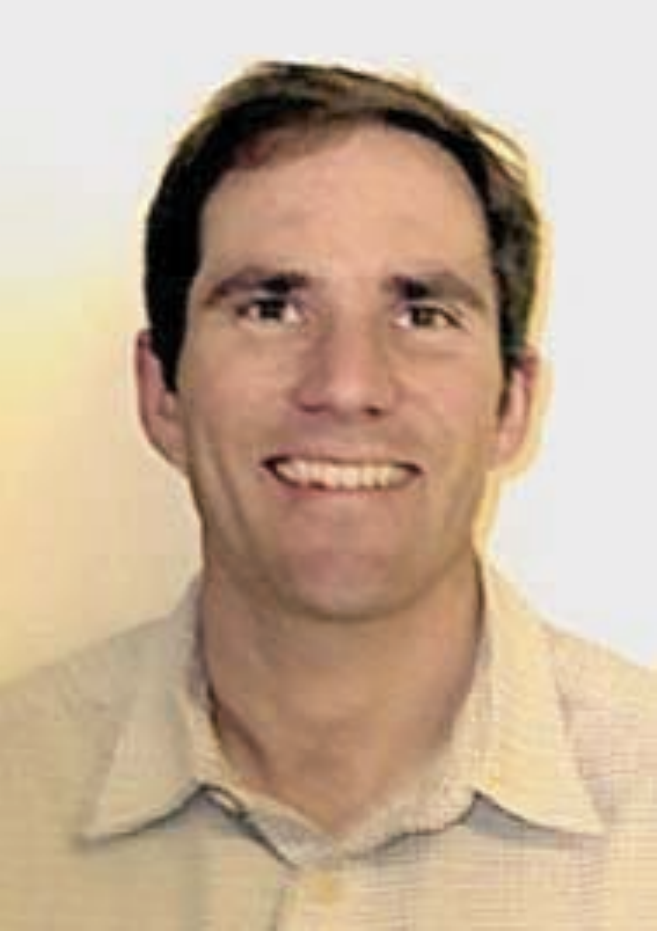}}]{John E. McInroy}

(M’91) received the B.S. degree
in electrical engineering from the University
of Wyoming, Laramie, in 1986, and both the
M.S. and Ph.D. degrees from Rensselaer Polytechnic
Institute, Troy, NY, in 1991. During his
graduate study, he was a Rensselaer Fellow, Xerox
Fellow, and NASA Graduate Student Researchers
Fellow.

He is currently a Professor of Electrical Engineering
at the University of Wyoming, Laramie,
USA. His interests are in the areas of computer vision, precision robotic
control, sensor fusion, machine intelligence, wind surfing, and
cross-country skiing.

Dr. McInroy is a Senior Member of IEEE.
\end{IEEEbiography}



\end{document}